\documentclass[a4paper]{article}

\usepackage{amsthm}

\newtheorem{thm}{Theorem}

\newtheorem{cor}{Corollary}
\newtheorem{Rem}{Remark}

\newtheorem{prop}{Proposition}[section]

\usepackage[toc,page]{appendix}
\usepackage{url}
\usepackage{pgfplots}
\usepackage{tikz}
\usepackage{amsmath}
\usepackage{mathtools}
\usepackage[a4paper, left=3cm, right=3cm, top=4cm]{geometry}
\usepackage{amssymb}
\usepackage{amsfonts}
\usepackage{scrlayer-scrpage}
\usepackage{accents}
\usepackage{ulem}
\clearscrheadfoot
\newcommand{\R}{\mathbb{R}}
\newcommand{\D}{\mathbb{D}}
\newcommand{\C}{\mathbb{C}}
\newcommand{\Z}{\mathbb{Z}}

\newcommand{\set}[1]{\left\{#1\right\}}

\newcommand{\abs}[1]{\left|#1\right|}

\newcommand{\norm}[1]{\left\|#1\right\|}

\newcommand{\ls}{\lesssim}

\def\Xint#1{\mathchoice
{\XXint\displaystyle\textstyle{#1}}
{\XXint\textstyle\scriptstyle{#1}}
{\XXint\scriptstyle\scriptscriptstyle{#1}}
{\XXint\scriptscriptstyle\scriptscriptstyle{#1}}
\!\int}
\def\XXint#1#2#3{{\setbox0=\hbox{$#1{#2#3}{\int}$ }
\vcenter{\hbox{$#2#3$ }}\kern-.6\wd0}}

\def\dashint{\Xint-}

\def\lf{\left}
\def\rg{\right}

\def\ds{\displaystyle}

\AtBeginDocument{
  \let\div\relax
  \DeclareMathOperator{\div}{div}
}

\usepackage{setspace}
\makeatletter
\newcommand{\MSonehalfspacing}{%
  \setstretch{1.44}
  \ifcase \@ptsize \relax 
    \setstretch {1.448}%
  \or 
    \setstretch {1.399}%
  \or 
    \setstretch {1.433}%
  \fi
}
\newcommand{\MSdoublespacing}{%
  \setstretch {1.92}
  \ifcase \@ptsize \relax 
    \setstretch {1.936}%
  \or 
    \setstretch {1.866}%
  \or 
    \setstretch {1.902}%
  \fi
}
\makeatother

\ohead{\pagemark}

\DeclareMathOperator{\sign}{sign}

\begin{document}

	\title{  Bergman-Bourgain-Brezis-type Inequality   }
	
	\author{ Francesca Da Lio,  Tristan Rivi\`ere, Jerome Wettstein\footnote{Department of Mathematics, ETH Zentrum, CH-8093 Z\"urich, Switzerland.}}
\maketitle
\begin{abstract}
In this note, we prove  a fractional version in $1$-D of the Bourgain-Brezis inequality \cite{bourgain1}. We show that such an inequality is equivalent to the fact that
a holomorphic function $f\colon\D\to\C$ belongs to the Bergman space ${\mathcal{A}}^2(\D)$, namely $f\in L^2(\D)$, if and only if 
$$\|f\|_{ L^1+ {H}^{-1/2}(S^1)}:=\limsup_{r\to 1^-}\|f(re^{i\theta})\|_{ L^1+ {H}^{-1/2}(S^1)}<+\infty.$$
Possible generalisations to the higher-dimensional torus are explored.
\end{abstract}
\medskip
{\noindent{\small{\bf Keywords.} Fractional Laplacian, holomorphic functions, Hardy spaces, Bergman spaces.}}\par

{\noindent{\small { \bf  MSC 2000.}  26A33, 32A10, 32A35, 
32A36 }}

\tableofcontents 
	
	\section{Introduction}
	
	In his pioneering work \cite{Riesz}, Riesz studied fine properties of the so-called Hardy spaces   ${\mathcal{H}}^p(\D)$,  which are the  spaces of holomorphic functions  \footnote{In 1915 Hardy observed that if $f$ is holomorphic in $\D$ then $r\mapsto M(r)=\int_0^{2\pi}|f(re^{i\theta})|^p d\theta$ is a nondcreasing function.
	 }  
	 $f\colon \D\to\C$ such that
	\begin{equation}\label{limcondr}
	\sup_{0<r<1}\int_0^{2\pi}|f(re^{i\theta})|^p d\theta<+\infty
	\end{equation}
	for $p>0$. Under condition \eqref{limcondr}, it is known that $f(e^{i\theta})$ exists and 
	\begin{equation}
	\lim_{r\to 1^-}\int_0^{2\pi}|f(re^{i\theta})-f(e^{i\theta})|^p d\theta=0~~\mbox{as well as}~~\lim_{r\to 1^-} f(re^{i\theta})=f(e^{i\theta}),
	\end{equation}
	for almost every $\theta.$ For $f\in {\mathcal{H}}^p(\D)$, one defines
	$\|f\|_{{\mathcal{H}}^p(\D)}:=\|f\|_{L^p(S^1)}.$
	\par
	We can independently  consider holomorphic functions in $L^2(\D)$ which corresponds to the well-known Bergman space ${\mathcal{A}}^2(\D)$
	\footnote{We recall that ${\mathcal{A}}^2(\D):=\{f\colon\D\to\C:~ f ~\mbox{holomorphic and $\|f\|_{L^2(\D)}<+\infty$}\}$}, see e.g \cite{DS}.\par
	The connection between Hardy spaces and the Bergman space ${\mathcal{A}}^2(\D)$ is given by the embedding ${\mathcal{H}}^1(\D) \hookrightarrow {\mathcal{A}}^2(\D)$ together with
	the estimate
	\begin{equation}\label{embHA}
	\|f\|_{L^2(\D)}\le C\|f\|_{{\mathcal{H}}^1(\D)}:=\|f\|_{L^1(S^1)}.
	\end{equation}
	In the case $\lim_{r\to 1^-}\|f(re^{i\theta})\|_{{H}^{-1/2}(S^1)}<+\infty$, then, by definition,  the following inequality holds   as well:
	\begin{equation}\label{embHA2}
	\|f\|_{L^2(\D)}\le C \|f\|_{H^{-1/2}(S^1)}:=\lim_{r\to 1^-}\|f(re^{i\theta})\|_{ {H}^{-1/2}(S^1)}.
	\end{equation}
	 In this note, we prove the following combination of \eqref{embHA} and \eqref{embHA2}:
	 
\begin{thm}\label{BBB}
Let $f\colon \D\to \C$ be an analytic function. Then $f$ belongs 
    to the Bergman space ${\mathcal{A}}^2(\D)$  if and only if  $$\|f\|_{ L^1+ {H}^{-1/2}(S^1)}:=\limsup_{r\to 1^-}\|f(re^{i\theta})\|_{ L^1+ {H}^{-1/2}(S^1)}<+\infty.$$
 Moreover, it holds
\begin{equation}\label{bergman}
\|f\|_{L^2(\D)}\le C\|f\|_{L^1+{H}^{-1/2}(S^1)}.\end{equation}
 \end{thm}
	 This type of inequalities takes its roots in the 
	  pioneering work  \cite{bourgain1}, where Bourgain and Brezis proved the following  striking result:
\begin{thm}[Lemma 1 in \cite{bourgain1}]\label{ThBB}
Let $u$  be  a $2\pi$-periodic function in $\R^n$  such that $\int_{\R^n} u\, dx=0$, and let 
$\nabla u=f + g,$  where $f\in \dot{W}^{-1,\frac{n}{n-1}}(\R^n)$\footnote{For $1< p < +\infty$, we will denote by $\dot{W}^{1,p}(\R^2)$ the homogeneous Sobolev space defined as the space of $f\in L^{1}_{loc}(\R^n)$ such that
$\nabla f\in L^{p}(\R^n)$ and by $\dot{W}^{-1,p^{\prime}}(\R^n)$ the corresponding dual space ($p^{\prime}$ is the conjugate of $p$). Every
function $f\in \dot{W}^{-1,p^{\prime}}(\R^n)$ can be represented as $f=\sum_{i=1}^n\partial_{x_i} f^j$ with $f^j\in L^{p'}(\R^n)$.} and $g\in L^1(\R^n)$  are  $2\pi$-periodic vector-valued functions. Then
\begin{equation}\label{BBineq}
\|u\|_{L^{\frac{n}{n-1}}}\le c \left(\|f\|_{\dot W^{-1,\frac{n}{n-1}}}+\|g\|_{L^1}\right).\
\end{equation}
 \end{thm}
 By duality, this implies the following corollary:
 \begin{cor}[Theorem 1 in  \cite{bourgain1}]\label{corBB}
For every $2\pi$-periodic function $h\in L^n(\R^n)$  with  $\int_{\R^n} h=0$, there exists a  $2\pi$-periodic $v\in \dot{W}^{1,n}\cap L^{\infty}(\R^n)$ 
  satisfying
$$\div v = h~~~\mbox{in $\R^n$}$$
and
\begin{equation}
	\|v\|_{L^{\infty}}+\|v\|_{\dot{W}^{1,n}}\le C(n)\|h\|_{L^n}.
\end{equation}
\end{cor}

One of the main result of this note is a fractional type Bourgain-Brezis inequality   on the circle  $S^1$.  More precisely,  we have the following:
 
\begin{thm}\label{nonlocalBB1D}
Let $u\in{\mathcal{D}}^{\prime}(S^1)$ be  such that $(-\Delta)^{\frac{1}{4}}u ,\mathcal{R}(-\Delta)^{\frac{1}{4}}u \in \dot{H}^{-\frac{1}{2}}(S^1) + L^{1}(S^1).$\footnote{We denote by ${\mathcal{R}}$ and  ${\mathcal{R}}_j$ the Riesz transform respectively on $S^1$ and  with respect to the $x_j$ variable on $T^n$, for  $j\in\{1,\ldots, n\}$ and  by  $\dot{H}^{-\frac{n}{2}}(T^n)$ the space of $f\in \mathcal{D}^\prime(S^1)$  such that $f=(-\Delta)^{n/4}g$, with $g\in L^2(T^n)$. Recall that $L^2_*(S^1):=\{u\in L^2(S^1):~~\dashint_{S^1}u =0\}$.}
Then  \\$u-\dashint_{S^1}u\in L^2_*(S^1)$ and the following estimate holds true:
\begin{equation}\label{fractineq}
\lf\|u-\dashint_{S^1}u\rg\|_{L^2}\le C\left(\|(-\Delta)^{1/4} u\|_{\dot{H}^{-1/2}(S^1)+L^1(S^1)}+\|{\mathcal{R}}(-\Delta)^{1/4} u\|_{\dot{H}^{-1/2}(S^1)+L^1(S^1)}\right),
\end{equation}
for some $C > 0$ independent of  $u$.
\end{thm}
We also show   the equivalence between Theorems  \ref{BBB} and \ref{nonlocalBB1D}, establishing the connection between fractional Bourgain-Brezis inequalities and Bergman spaces.\\

\par
We would like to add some comments  about 	 Bourgain-Brezis' inequality. 
 Bourgain-Brezis' inequality  in its general form is of interest in the study of the PDE $\operatorname{div} Y = f$ for $f \in L^{n}_{\ast}(T^{n})$, where for finite $p\ge 1$, $L^{p}_{\ast}(T^n)$ denotes the Banach subspace of $L^p$-functions with vanishing mean over the torus. Precisely, they found that $Y$ can be chosen to be continuous and in $\dot{W}^{1,n}(T^n)$, a result which is non-trivial due to the fact that $\dot{W}^{1,n}(T^n)$ does not continuously embed into $L^{\infty}(T^n)$. 
 The key ingredient in the proof is a duality argument based on an estimate similar to \eqref{BBineq} and some general results from functional analysis regarding closedness properties of the image space. This motivates the general interest in inequalities of the same type, as improved regularity results in limit cases can be invaluable. Indeed, later, such estimates have been considered by a variety of further authors, see \cite{bourgain2} again by Bourgain and Brezis for results regarding Hodge decompositions, \cite{mazya} due to Maz'ya for an inequality on $H^{1-\frac{n}{2}}(\mathbb{R}^{n})$ leading to a different existence result for the PDE $\operatorname{div} Y = f$ and \cite{mironescu} by Mironescu using a PDE-approach to a $2$D-special case involving explicitly computing fundamental solutions with appropriate boundedness properties.\par
	In \cite{dalioriv},  the first two authors provide an alternative proof of  \eqref{BBineq} in dimension $2$ without 
	 assuming the periodicity of the function $u$. The proof is related to some compensation phenomena observed first in \cite{De} in the analysis of $2$-dimensional perfect incompressible fluids and then  also  applied  by the second author in  \cite{Riv3} in  the analysis of {\it isothermic surfaces}. 
	 For an overview of the results in the literature regarding variations of  Theorem \ref{ThBB} and Corollary \ref{corBB}, we refer for instance to \cite{vs}.
	  The inequality  \eqref{BBineq} also represents the first key step in the study of the regularity of  $ L^2(\D,{\R}^n)$ solutions
	 $u$  to a linear elliptic system of the following form
\begin{eqnarray}
\label{intro-6}
\mbox{div}\,(S\,\nabla u)&=&\sum_{j=1}^n\mbox{div}\,(S_{ij}\,\nabla u^j)=\sum_{j=1}^n\sum_{\alpha=1}^2\frac{\partial}{\partial {x_{\alpha}}}\left(S_{ij}u_{x_{\alpha}}^j\right)=0,
\end{eqnarray}
where $S$ is a $ W^{1,2}(\D)$ symmetric $n\times n$ matrix, such that $S^2=id_n$, as seen in \cite{dalioriv}.\\

We would like to mention some results on Riesz potentials showing that 
  the $1$-dimensional case plays a particular role in the $L^1$-estimates for Riesz potentials. More precisely, one can deduce from the results in \cite{sw} that for all
  $0<\alpha<1$,  we have:
  \begin{equation}\label{steinineq}\|I_{\alpha}u\|_{L^{\frac{1}{1-\alpha}}}\le C(\|{\mathcal{R}}u\|_{L^1}+\|u\|_{L^1}),\end{equation}
  for all $u$ in the Hardy space ${\mathcal{H}}^1(\R)$. 
  It follows in particular that:
$$\| u\|_{L^{\frac{1}{1-\alpha}}}\le C(\|{\mathcal{R}}(-\Delta)^{\alpha/2}u\|_{L^1}+\|(-\Delta)^{\alpha/2}u\|_{L^1}).$$
In \cite{ssv}, the authors show that if $N\ge 2$ and $0<\alpha<N$,  then there is a constant $C=C(\alpha,N)>0$ such that
\begin{equation}\label{estschr}
\|I_{\alpha}u\|_{L^{\frac{N}{N-\alpha}}}\le C\|{\mathcal{R}}u\|_{L^1} \end{equation}
for all $u\in C^{\infty}_c(\R^N)$, such that ${\mathcal{R}}u\in L^1(\R^N).$
The estimate \eqref{estschr} is however false  in $1$-D, as seen in \cite{ssv}.
\par
The inequality \eqref{fractineq} generalizes  the inequality \eqref{steinineq} in the case  $\alpha=1/2$ and the counter-example
  in \cite{ssv} for the estimate \eqref{estschr} in $1$-D shows that the estimate \eqref{fractineq}  is in some sense optimal.\\
  
  The paper is organized as follows: In section \ref{prel}, we recall the definitions of the fractional    Laplacian on the unit circle and on the torus and of Clifford Algebras. In section \ref{proof}, we provide two distinct proofs of Theorem \ref{nonlocalBB1D} .
In section \ref{equiva},  we establish the equivalence of    Theorem \ref{BBB}  and Theorem \ref{nonlocalBB1D}.
 In section \ref{torus},  we extend the fractional Bourgain-Brezis inequality using Clifford algebras   to the $n$-dimensional torus $T^n$. 
In section \ref{exist}, we prove 
  existence results for certain fractional PDE-operators in the same  spirit as Corollary \ref{corBB}. Lastly, in section \ref{appendix}, we provide a proof of the inequalities \eqref{embHA} and \eqref{embHA2}.

  \par

	\section{Preliminaries}\label{prel}
	\subsection{Fractional    Laplacian on the unit circle and on the torus}
	Before we enter the discussion and the proofs of the main results, let us recall a few notions essential in our later arguments. We mainly focus on fractional Laplacians, fractional Sobolev spaces and Clifford algebras.\par
	Throughout this note, we shall denote by $T^n$ the torus of dimension $n \in \mathbb{N}$. This means:
	\begin{equation}
		T^n = \underbrace{S^1 \times \ldots \times S^1}_\text{$n$ times} = \R^n / (2\pi \Z)^n
	\end{equation}
	where  $  S^1=\R/2\pi\Z$.
	We denote by $\mathcal{D}( T^n):=C^\infty(T^n)$ the Fr\'echet space of smooth functions on $T^n$ and by $\mathcal{D}'(T^n)$ its topological dual. The natural duality paring is denoted by $\langle \cdot, \cdot \rangle$.
 
For $u\in\mathcal{D}'(T^n)$ and $m\in\mathbb{Z}^{n}$, we define the Fourier coefficients of $u$ as follows:
\begin{equation}
		\widehat{u}(m) := \frac{1}{(2 \pi)^n} \int_{T^n} u(x) e^{-i \langle m, x \rangle} dx = \Big{\langle} u, e^{-i \langle m, \cdot \rangle} \Big{\rangle}.
	\end{equation}  
The Fourier coefficients completely determine $u$ as a distribution on $T^n$ and convergence in the sense of distributions obviously implies convergence of the Fourier coefficients. Notice that, for all $u\in\mathcal{D}'(T^n)$, there exists some $N>0$ such that $\abs{\hat u(m)}\ls(1+\abs{m})^N$. Moreover, we recall that $v\in C^\infty(T^n)$ if and only if the Fourier coefficients $\hat v(m)$ have rapid decay, i.e. $\sup_m(1+\abs{m})^N\abs{\hat v(m)}<\infty$ for all $N>0$.

Given $s\in\R$, we define the non-homogeneous and homogeneous Sobolev spaces respectively by
\[ H^s( T^n):=\set{v\in\mathcal{D}'(T^n):\norm{v}_{H^s}^2:=\sum_{k\in\mathbb{Z}^{n}}(1+\abs{k}^2)^s\abs{\hat v(k)}^2<\infty}, \]
and 
\[ \dot{H}^s( T^n):=\set{v\in\mathcal{D}'(T^n):\norm{v}_{\dot H^s}^2:=\sum_{k\in\mathbb{Z}^{n}}|k|^{2s}\abs{\hat v(k)}^2<\infty}, \]
where  $\mathcal{D}'(T^n)$ is again the space of distributions on $T^n$. Notice that if $s = 0$, we have $L^{2}(T^n) = H^{0}(T^n)$ and $L^{2}_{\ast}(T^n) \simeq \dot H^{0}(T^n)$.

	An important family of operators throughout our considerations are the so-called \textit{fractional Laplacians}. Let $s > 0$ be real, then we define for $u: T^n \to \mathbb{C}$ smooth the $s$-Laplacian of $u$ by the following multiplier property:

	\begin{equation}
		\widehat{(-\Delta)^{s} u}(\xi) = \sum_{m\in\mathbb{Z}^{n}}|m|^{2 s} \widehat{u}(m)e^{i\langle m,\xi\rangle} , \quad \forall \xi \in T^n.
	\end{equation}
	Clearly, this definition can immediately be extended to the spaces $H^{s}(T^{n})$ or even $\mathcal{D}^\prime (T^n)$ as a multiplier operator on Fourier coefficients, possibly defining merely a distribution on $T^n$.
	Finally, we recall the definition of the $j$-Riesz transform on $T^n$ as a multiplier operator:
	
	\begin{equation}
		\widehat{{\mathcal{R}}_ju}(\xi) = \sum_{m\in\mathbb{Z}^{n}}i \frac{m_j}{|m|} \widehat{u}(m)e^{i\langle m,\xi\rangle} , \quad \forall \xi \in T^n.
	\end{equation}
	In particular, in the case $n = 1$, we have:
	\begin{equation}
		\widehat{{\mathcal{R}}u}(\xi) = \sum_{m\in\mathbb{Z}}i \sign(m) \widehat{u}(m)e^{i m \cdot \xi} , \quad \forall \xi \in S^1.
	\end{equation}
	
	\subsection{Clifford Algebras}
	The material covered here is due to \cite{gilbert} and \cite{hamilton} and we refer to them for further details on the topics introduced. For the remainder of this subsection, let $\mathbb{K} \in \{ \mathbb{R}, \mathbb{C} \}$ denote a scalar field and $V$ a finite dimensional $\mathbb{K}$-vector space. Let $Q: V \to \mathbb{K}$ be a map, such that:
	\begin{itemize}
		\item[1.)] For all $\lambda \in \mathbb{K}$ and $v \in V$, we have: $Q(\lambda v ) = \lambda^2 \cdot Q(v)$.
		\item[2.)] The map $B(v,w) := \frac{1}{2} \big{(} Q(v+w) - Q(v) - Q(w) \big{)}$ defines a $\mathbb{K}$-bilinear map on $V \times V$.
	\end{itemize}
	Such a $Q$ will be called a \textit{quadratic form} and the pair $(V,Q)$ a \textit{quadratic space}. Standard examples include real vector spaces equipped with scalar products, but not complex vector spaces with scalar products due to complex anti-linearity in the second argument. Inspired by this example, we say that a basis $e_1, \ldots, e_n$ of a quadratic space $(V,Q)$ is \textit{$B$-orthonormal}, if for all $j \in \{ 1, \ldots, n \}$, we have $| Q(e_j) | = 1$ as well as:
	\begin{equation}
		B(e_j, e_k) = 0, \quad \forall j \neq k \in \{ 1, \ldots, n \}.
	\end{equation}
	Given such a quadratic space $(V, Q)$, we call a pair $(\mathcal{A}, \nu)$ a \textit{Clifford algebra} for $(V,Q)$, if the following holds, see \cite[p.8, (2.1)]{gilbert}:
	\begin{itemize}
		\item[i.)] $\mathcal{A}$ is an associative algebra with unit $1$ and $\nu: V \to \mathcal{A}$ is $\mathbb{K}$-linear and injective.
		\item[ii.)] $\mathcal{A}$ is generated as an algebra by $\nu(V)$ and $\mathbb{K} \cdot 1$.
		\item[iii.)] For every $v \in V$, we have: $\nu(v)^{2} = - Q(v) \cdot 1$
	\end{itemize}
	An important immediate corollary of the definition is the following commutation relation:
	\begin{equation}
		\nu(v) \nu(w) + \nu(w) \nu(v) = - 2 B(v,w) \cdot 1, \quad \forall v, w \in V.
	\end{equation}
	Thus, pairs of {\em orthogonal} vectors with respect to $B$ anti-commute as elements in $\mathcal{A}$. We usually omit explicitly mentioning $\nu$ and therefore identify $v$ with $\nu(v)$, which is justified due to $\nu$ being injective.\\
	
	For the remainder of the section, let us focus on $(V,Q)$ non-degenerate, i.e. for all $v \in V$, there is a $w \in V$, such that $B(v,w) \neq 0$. In this case, there actually exists a basis $e_1, \ldots, e_n$, where $n = \operatorname{dim}_{\mathbb{K}} V$, orthonormal with respect to $B$ and, consequently, such that:
	\begin{equation}
		e_j e_k + e_k e_j = \pm 2 \delta_{jk} \cdot 1, \quad \forall j,k \in \{ 1, \ldots, n\},
	\end{equation}
	(see e.g. Theorem 1.5 in \cite{gilbert}). The signs are determined by the signature of the quadratic form $Q$ and may vary for different choices $j,k$. Provided $\mathbb{K} = \mathbb{C}$, we may assume that all signs are the same, see \cite{gilbert}.\par
	It can be shown that every Clifford algebra has $\mathbb{K}$-dimension at most $2^{n}$. If the dimension is equal to $2^n$, the Clifford algebra is called \textit{universal}.\footnote{This definition is justified, as universal Clifford algebras $\mathcal{A}$ have an extension property for linear maps from $V$ to any Clifford algebra respecting the characteristic multiplication relation in $\mathcal{A}$, see \cite{gilbert}.} An important result in \cite[Thm. 2.7]{gilbert} states that there always exists a universal Clifford algebra for any given quadratic space. Moreover, there exist explicit descriptions of all universal Clifford algebras up to isomorphisms in terms of matrices, see \cite{gilbert}.\\
	
	To conclude this brief treatment of Clifford algebras, let us provide an explicit example: Let $V = \mathbb{C}^{n}$, $\mathbb{K} =\mathbb{C}$ and define $Q$ as follows:
	\begin{equation}
		Q(z_1, \ldots, z_n) := \sum_{j=1}^{n} z_j^2, \quad \forall (z_{1}, \ldots, z_n) \in \mathbb{C}^{n}.
	\end{equation}
	It is clear that $(V,Q)$ is a non-degenerate quadratic space, as $B$ is the standard scalar product up to a complex conjugation in the second argument. In this case, the standard basis $e_1, \ldots, e_n$ already is $B$-orthonormal. Thus, we have:
	\begin{equation}
		e_j e_k + e_k e_j = - 2 \delta_{jk} \cdot 1, \quad \forall j,k \in \{ 1, \ldots, n \}.
	\end{equation}
	The universal Clifford algebra is then spanned by the finite products $e_{\alpha}$ of the basis elements, where $\alpha \subset \{ 1, \ldots, n \}$ is an ordered subset and we define:
	$$e_\alpha = \prod_{j \in \alpha} e_{j}$$
	In particular, $e_{\emptyset} = 1$ by definition. It can be seen that every complex universal Clifford algebra associated with a non-degenerate quadratic space of dimension $n$ is isomorphic to this one, see \cite{gilbert} and the definition of universal Clifford algebra presented there.\\
	
	Lastly, let us introduce a few definitions from Chapter 1, Section 7 in \cite{gilbert}: We may identify the universal Clifford algebra ${\mathcal{A}}$ as a vector space with $\mathbb{K}^{2^{n}}$, if $\operatorname{dim}_{\mathbb{K}} V = n$. This allows us to generalize the natural scalar product-induced norm on $\mathbb{K}^{2^{n}}$ to the Clifford algebra and we shall denote this norm by $\| \cdot \|$. Moreover, there is a notion of conjugation on Clifford algebras defined by:
	\begin{equation}
	\label{cliffconjdef}
		\overline{e_{j_1} \ldots e_{j_{k}}} := (-1)^{k} Q(e_{j_1}) \ldots Q(e_{j_k}) \cdot e_{j_k} \ldots e_{j_1} = (-1)^{\frac{k(k+1)}{2}} Q(e_{j_1}) \ldots Q(e_{j_k}) \cdot e_{j_1} \ldots e_{j_k},
	\end{equation}
	and extending linearily. If $\mathbb{K} = \mathbb{C}$, we also conjugate the complex coefficients in the usual manner, i.e. we extend complex anti-linearily. We highlight the following key property of the conjugation:
	\begin{equation}
		\overline{xy} = \overline{y} \cdot \overline{x}, \quad \forall x,y \in \mathcal{A}.
	\end{equation}
	This is  due to the inversion of factors in  \eqref{cliffconjdef}. We emphasise that the definition in \eqref{cliffconjdef} is precisely made with the identity below in mind:
	\begin{equation}
		\overline{e_{j_1} \ldots e_{j_{k}}} \cdot e_{j_1} \ldots e_{j_{k}} = 1.
	\end{equation}
	
	The following property will be useful later as well: Let $x \in \mathcal{A}$ be given and denote by $P_{0}$ the linear projection of an element in the Clifford algebra to the coefficient associated with the neutral element $1$. More precisely, $P_0: \mathcal{A} \to \mathbb{K}$ is the following linear map:
	$$P_0 \Big{(} \sum_{\alpha} x_{\alpha} e_{\alpha} \Big{)} = x_{\emptyset}$$
	We have by a direct computation:
	\begin{align}
		P_{0}(\overline{x}x)	&= \sum_{\alpha \subset \{1, \ldots, n\}} \overline{x_{\alpha}} x_{\alpha} \notag \\
						&= \| x \|^2,
	\end{align}
	where we wrote explicitly $x = \sum_{\alpha \subset \{1, \ldots, n\}} x_{\alpha} e_{\alpha}$ with $x_{\alpha} \in \mathbb{K}$. It suffices to observe that $e_\alpha \cdot e_\beta$ has non-vanishing contribution in the $e_{\emptyset} = 1$-direction, if and only if $\alpha = \beta$. The formula then follows.
 
	\section{Fractional Bourgain-Brezis inequality on the unit circle ${\mathbf{S}^1}$}\label{proof}
	In this section, we provide two distinct proofs of Theorem \ref{nonlocalBB1D}. The first proof is in the spirit of the one presented in \cite{bourgain1}, while the second one is inspired by that in \cite{dalioriv} and is based on some particular compensation phenomena. We assume for simplicity that $u$ is real valued (the proof for  complex-valued function is completely analogous, see Remark \ref{complex}).\par
	First, we would like to observe  that if $u\in C^{\infty} (S^1)$, then by definition: 
\begin{eqnarray}\label{nlest}
\ds\lf\|u-\dashint_{S^1}u\rg\|_{L^2}&\le &C \|(-\Delta)^{1/4} u\|_{\dot{H}^{-1/2}(S^1)}.
\end{eqnarray}
On the other hand, as we have already observed in the introduction,  we also have\footnote{Actually, an even sharper inequality than \eqref{nlestbis} holds true with $L^2(S^1)$ being replaced by the smaller Lorentz space $L^{2,1}(S^1)$.}:
\begin{eqnarray}
\label{nlestbis}
\ds\lf\|u-\dashint_{S^1}u\rg\|_{L^2(S^1)}&\le &C\left(\|(-\Delta)^{1/4} u\|_{L^1(S^1)}+\|{\mathcal{R}}(-\Delta)^{1/4} u\|_{L^1(S^1)}\right)\simeq \|(-\Delta)^{1/4} u\|_{{\mathcal H}^1(S^1)}.
\end{eqnarray}
  	 
	 \subsection{A first proof of Theorem \ref{nonlocalBB1D}}
  
Let us suppose that  $u\in C^{\infty} (S^1)$,  $\dashint_{S^1}u=0$. We assume for simplicity that  $u$ is real-valued, (see Remark \ref{complex} for the complex-valued case).  The proof below follows the main arguments of the original proof by Bourgain and Brezis.
We write:
\begin{equation}\label{estnl1}
\left\{\begin{array}{c}
 (-\Delta)^{1/4} u=f^{1}+g^{1}\\[5mm]
 {\mathcal{R}} (-\Delta)^{1/4} u=f^{2}+g^{2}
 \end{array}\right.
 \end{equation}
 where $f^{1},f^{2}\in \dot{H}^{-1/2}(S^1)$, $g^{1},g^{2}\in L^1(S^1)$.
 We set $u=\sum_{n\in\Z^*} u_n e^{in \theta}$. Since $u$ is real-valued, it holds $\bar u_n=u_{-n}.$
 We see:
 \begin{eqnarray}\label{estnl2}
 \sum_{n\in\Z^*} |u_n|^2&=& \sum_{n\in\Z^*} |n|^{1/2}u_n\frac{u_{-n}}{|n|^{1/2}}= \sum_{n\in\Z^*} \frac{f^{1}_n+g^{1}_n}{|n|^{1/2}}u_{-n},
 \end{eqnarray}
  \begin{eqnarray}\label{estnl3}
 \sum_{n\in\Z^*}\frac{f^{1}_n\, u_{-n}}{|n|^{1/2}}&\le & \left[ \sum_{n\in\Z^*}\frac{|f^{1}_n|^2}{|n|}\right]^{1/2}\left[ \sum_{n\in\Z^*}|u_n|^2\right]^{1/2},
 \end{eqnarray}
 \begin{eqnarray}\label{estnl4}
 \sum_{n\in\Z^*}\frac{g^{1}_n\, u_{-n}}{|n|^{1/2}}&= & \sum_{n>0}\frac{g^{1}_n\, u_{-n}}{|n|^{1/2}}+ \sum_{n<0}\frac{g^{1}_n\, u_{-n}}{|n|^{1/2}}.
 \end{eqnarray}
 Observe that by definition of the Riesz transform:
  \begin{eqnarray}\label{estnl5}
  {\mathcal{R}} (-\Delta)^{1/4} u&=&i\left[- \sum_{n<0} |n|^{1/2} u_n e^{in\theta}+\sum_{n>0} |n|^{1/2} u_n e^{in\theta}\right].
  \end{eqnarray}
  Therefore:
  \begin{equation}\label{estnl6}
  u_n = \left\{\begin{array}{cc}
  \frac{f_n^2+g_n^2}{-i|n|^{1/2}}&~~~\mbox{if $n<0$}\\[5mm]
  \frac{f_n^2+g_n^2}{i|n|^{1/2}}&~~~\mbox{if $n>0$}
  \end{array}\right.\end{equation}
  By combining \eqref{estnl4} and \eqref{estnl6}, we obtain:
  \begin{eqnarray}\label{estnl7}
 \sum_{n\in\Z^*}\frac{g^{1}_n\, u_{-n}}{|n|^{1/2}}&= &\sum_{n>0}g_n^1\ \frac{f_{-n}^2+g_{-n}^2}{-i|n|}+\sum_{n<0}g_n^1\ \frac{f_{-n}^2+g_{-n}^2}{i|n|}.
 \end{eqnarray}
 Let us estimate the different parts of the sum \eqref{estnl7} individually:\\
 
 \noindent
{ 1.  We first estimate}
 \begin{eqnarray}\label{estnl8}
  \sum_{n\in\Z^*}{\mbox{sign$(n)$}}\frac{g^{1}_n\ f^2_{-n}}{|n|}&=& \sum_{n\in\Z^*}{\mbox{sign$(n)$}}\frac{| n |^{1/2}u_n-f^1_n}{|n|^{1/2}}\frac{f_{-n}^2}{|n|^{1/2}}\nonumber\\[5mm]
  &\le& \left( \sum_{n\in\Z^*}|u_n|^2\right)^{1/2}  \left( \sum_{n\in\Z^*}\frac{|f^2_n|^2}{|n|}\right)^{1/2} +\left( \sum_{n\in\Z^*}\frac{|f^1_n|^2}{|n|}\right)^{1/2}\left( \sum_{n\in\Z^*}\frac{|f^2_n|^2}{|n|}\right)^{1/2}\nonumber\\[5mm]
  &\le& \|u\|_{L^2}\|f^2\|_{\dot{H}^{-1/2}}+ \|f^1\|_{\dot{H}^{-1/2}}\|f^2\|_{\dot{H}^{-1/2}}.
  \end{eqnarray}
{ 2. It remains to estimate}
 $$ \sum_{n\in\Z^*} {\mbox{sign$(n)$}}\frac{g^{1}_n\, g^2_{-n}}{i|n|}.$$
 For this purpose, we consider the following operator:
 \begin{eqnarray*}
 {\mathbf{A}}\colon L^1(S^1)\times L^1(S^1)&\to &\C ,
 ~~~~~~(g^1,g^2)\mapsto  \sum_{n\in\Z^*}{\mbox{sign$(n)$}}\frac{g^{1}_n\, g^2_{-n}}{i|n|}.
 \end{eqnarray*}
 {\bf Claim 1. } The operator $ {\mathbf{A}}$ is continuous, i.e. we have the following estimate:
 \begin{equation}\label{estn9} | {\mathbf{A}}(g^1,g^2)|\le C \|g^1\|_{L^1}\|g^2\|_{L^1}.\end{equation}
 {\bf Proof of Claim 1.} It is sufficient to prove the claim in the case where $g^1$ and $g^2$ are arbitrary Dirac-delta measures.\footnote{We recall that the  linear span of Dirac measures is dense in the space of Radon measures ${\mathcal{M}}(S^1)$ equipped with the weak-* topology.} Therefore, we consider $g^{1}=\sum_{i\in I}\lambda_i\delta_{a_i}$ and  $g^{2}=\sum_{j\in J}\mu_j\delta_{b_j}$. We have 
 $\|g^{1}\|_{{\mathcal{M}}(S^1)}=\sum_{i\in  I}|\lambda_i|,$  $\|g^{2}\|_{{\mathcal{M}}(S^1)}=\sum_{j\in  J}|\mu_j|.$ 
 By bilinearity, we deduce:
 \begin{eqnarray}\label{estnl10}
| {\mathbf{A}}(g^1,g^2)|&=&| {\mathbf{A}}(\sum_{i\in I}\lambda_i\delta_{a_i},\sum_{j\in J}\mu_j\delta_{b_j})|\nonumber\\[5mm]
&\le& \sum_{i\in I, j\in J}|\lambda_i||\mu_j|| {\mathbf{A}}(\delta_{a_i},\delta_{b_j})|\nonumber\\[5mm]
&\le& \sup_{(a,b)\in S^1\times S^1}| {\mathbf{A}}(\delta_a,\delta_{b})| \sum_{i\in I}|\lambda_i| \sum_{j\in J}|\mu_j|\nonumber\\[5mm]
&=& \sup_{(a,b)\in S^1\times S^1}| {\mathbf{A}}(\delta_{a},\delta_{b})|\|g^1\|_{{\mathcal{M}}(S^1)}|\|g^2\|_{{\mathcal{M}}(S^1)}.
\end{eqnarray}
If $ \sup_{(a,b)\in S^1\times S^1}| {\mathbf{A}}(\delta_{a},\delta_{b})|<+\infty$, then the claim holds  for  linear combinations of Dirac measures.
By a density argument, we get the claim 1 for arbitrary $g^1,g^2\in L^1(S^1).$ Hence, claim 1 is a consequence of the following:

\noindent{\bf Claim 2}.  $ \sup_{(a,b)\in S^1\times S^1}| {\mathbf{A}}(\delta_{a},\delta_{b})|<+\infty$.\par
\noindent{\bf Proof of Claim 2.}  For $g^1=\delta_a$ and $g^2=\delta_b$ , we have
$g_n^1=e^{ina}$ and $g^2_n=e^{inb}$. 
In this case, we observe:
\begin{eqnarray}\label{estnl111}
{\mathbf{A}}(\delta_a,\delta_b)&=&\sum_{n\in Z^*}\mbox{sign$(n)$}\frac{g_n^1 g_{-n}^2}{i|n|}\nonumber\\
&=&\sum_{n\in \Z^*}\mbox{sign$(n)$}\frac{e^{in(a-b)}}{{i|n|}}=2 \sum_{n>0}\frac{\sin(n(a-b))}{n }<+\infty.\footnotemark
\end{eqnarray}
This proves claim 2 and from \eqref{estnl111}, we can  deduce claim 1 as well.
\footnotetext{The value of such a series is deduced from the Fourier series of $f(x)=\frac{x}{2\pi}$ for $0<x<2\pi$ and $f(x+2\pi)=f(x)$.}\\

\noindent
By combining \eqref{estnl2}-\eqref{estnl10}  we get
\begin{eqnarray}\label{estnl11}
\|u\|^2_{L^2}&\lesssim & \|u\|_{L^2}\left(\|f^{1}\|_{\dot{H}^{-1/2}}+\|f^{2}\|_{\dot{H}^{-1/2}}\right)+ \|f^{1}\|_{\dot{H}^{-1/2}}\|f^{2}\|_{\dot{H}^{-1/2}}+C\|g^1\|_{L^1}\|g^2\|_{L^1} \nonumber\\[5mm]
&\lesssim & \frac{1}{2} \|u\|_{L^2}^2+\frac{1}{2} \left(\|f^{1}\|^2_{\dot{H}^{-1/2}}+\|f^{2}\|^2_{\dot{H}^{-1/2}}\right)+ \|f^{1}\|_{\dot{H}^{-1/2}}\|f^{2}\|_{\dot{H}^{-1/2}}+C\|g^1\|_{L^1}\|g^2\|_{L^1}  \nonumber \\[5mm]
&\lesssim&\frac{1}{2} \|u\|_{L^2}^2+\left(\|f^{1}\|^2_{\dot{H}^{-1/2}}+\|f^{2}\|^2_{\dot{H}^{-1/2}}\right)+\frac{1}{2}\left(\|g^1\|^2_{L^1}+\|g^2\|^2_{L^1}\right).
\end{eqnarray}
This estimate permits us to conclude the proof of Theorem \ref{nonlocalBB1D}. 
Since $f^1,f^2,g^1,g^2$ were arbitrary, one can deduce \eqref{fractineq}.
In the general case where
$u\in{\mathcal{D}}^{\prime}(S^1)$, one argues by approximation (see section \ref{secproof} for further details).~~\hfill$\Box$

\bigskip
 	\medskip
 	
	\subsection{A second proof of Theorem \ref{nonlocalBB1D}}\label{secproof}

	As in the first proof, we will show the following:
		Let $u \in \mathcal{D}'(S^{1})$ be  such that:
		\begin{align}
		\label{condition01}
			(-\Delta)^{\frac{1}{4}}u 			&= f_{1} + g_{1} \\
		\label{condition02}
			\mathcal{R}(-\Delta)^{\frac{1}{4}}u 	&= f_{2} + g_{2},
		\end{align}
		where $f_{1}, f_{2} \in \dot{H}^{-\frac{1}{2}}(S^{1})$ and $g_{1}, g_{2} \in L^{1}(S^{1})$. Under these conditions, we prove:
		\begin{equation}
			u - \int_{S^{1}} u dx \in L^{2}_{\ast}(S^{1})=\lf\{u\in L^2(S^1):~~\dashint_{S^1} u=0\rg\},
		\end{equation}
		together with the following estimate:
		\begin{equation}
		\label{estimateforuinthm}
			\Big{\|} u - \int_{S^{1}} u dx \Big{\|}_{L^{2}} \leq C \big{(} \| f_{1} \|_{\dot{H}^{-\frac{1}{2}}} + \| f_{2}\|_{\dot{H}^{-\frac{1}{2}}} + \| g_{1}\|_{L^{1}} + \| g_{2}\|_{L^{1}} \big{)},
		\end{equation}
		where $C > 0$ is independent of $f_{1}, f_{2}, g_{1}, g_{2}$ and $u$.
We may assume for simplicity that $u$ is real-valued (see Remark \ref{complex} for the complex-valued case).\par

		Firstly, observe that it suffices to consider the case:
		\begin{equation}
			\int_{S^{1}}  udx = 2\pi \cdot \hat{u}(0) = 0,
		\end{equation}
		by merely changing $u$ by a constant.  Similarly, by the conditions in \eqref{condition01} and \eqref{condition02}, we see that $f_{j}, g_{j}$ have vanishing integral over $S^{1}$ and consequently vanishing Fourier coefficient for $n = 0$.\footnote{It would be possible to treat $f_j, g_j$ with non-vanishing integral, i.e. treat the case $(-\Delta)^{{1/4}} u, \mathcal{R} (-\Delta)^{{1/4}} u \in L^{1} + H^{-{1/2}}(S^1)$ by reducing to vanishing Fourier coefficient at $n = 0$: We have by the conditions    $\widehat{f_j}(0) =-\widehat{g_j}(0)$. Note that $| \widehat{g_j}(0) | \lesssim \| g_j \|_{L^1}$. Note that $\| f_j \|_{{H}^{-1/2}}^2 \simeq | \widehat{f_{j}}(0) |^2 + \| \tilde{f}_j \|_{\dot{H}^{-1/2}}^2 $, where $\tilde{f}_j$ denotes the corrected $f_j$ with vanishing $0$th Fourier coefficient. Thus, we could reduce to the case of vanishing integral.}
		 For now, let us assume that $u, f_{j}, g_{j}$ are all smooth on $S^{1}$. The general case can be dealt with using convolution with an appropriate smoothing kernel and approximation arguments as specified at the end of the proof.\\
		
		First, let us define the following operators on $\mathcal{D}'(S^{1})$:
		\begin{align}
		\label{operatord}
			Dv 			&:= (-\Delta)^{\frac{1}{4}} \big{(} Id + \mathcal{R} \big{)}v \\
		\label{operatordbar}
			\overline{D}v 	&:= (-\Delta)^{\frac{1}{4}} \big{(} Id - \mathcal{R} \big{)}v,
		\end{align}
		for every $v \in \mathcal{D}'(S^{1})$. Consequently, using \eqref{condition01} and \eqref{condition02}, we have:
		\begin{align}
		\label{rearr1}
			Du 			&= f_{1} + f_{2} + g_{1} + g_{2} = f + g \\
		\label{rearr2}
			\overline{D}u 	&= f_{1} - f_{2} + g_{1} - g_{2} = \tilde{f} + \tilde{g}.
		\end{align}
		Let us calculate the Fourier multipliers associated with $D, \overline{D}$. For every $n \in \mathbb{Z}$, we have:
		\begin{align}
			\mathcal{F}\big{(} Dv\big{)}(n)			&= |n|^{\frac{1}{2}}(1+ i \sign(n))\hat{v}(n) \\
			\mathcal{F}\big{(} \overline{D}v\big{)}(n)	&= |n|^{\frac{1}{2}}(1- i \sign(n))\hat{v}(n)
		\end{align}
{\bf Claim 1:} Given $f\in \dot{H}^{-1/2}(S^1)$, there is 	  a real-valued function $F \in L^{2}_{\ast}(S^{1})$\footnote{The observation that $F$ may be chosen real-valued is due to $\widehat{F}(-n) = \overline{\widehat{F}(n)}$ for all $n$.}, such that $DF = f$.
\par
\noindent
{\bf Proof of the Claim 1}  In order to solve $DF=f$,   we should have:
	\begin{equation}
	\label{fourierF}
		\hat{F}(n) = \frac{1}{1 + i \sign(n)} \frac{\hat{f}(n)}{\sqrt{|n|}}, \quad \text{if } n \neq 0 \\
	\end{equation}
	Using the fact that the $L^{2}$-norm of $F$ can be characterized in terms of the $l^{2}$-norm of the Fourier coefficients, we obtain:
	\begin{align}
	\label{fouriercoeffest}
		\| F \|_{L^{2}}^{2}	&= \sum_{n \neq 0} \frac{1}{| 1 + i \sign(n) |^{2}} \frac{| \hat{f}(n) |^{2}}{|n|} \notag \\
						&\leq \sum_{n \neq 0} \frac{| \hat{f}(n) |^{2}}{|n|} \notag \\
						 						&=  \| f \|_{\dot{H}^{-\frac{1}{2}}}^{2},
	\end{align}

	\noindent
	where we used the definition of the $\dot{H}^{-\frac{1}{2}}$-norm. Observe that a converse inequality could be obtained along the same lines.
	
	Next, by defining $\tilde{u} := u - F$, we observe that due to \eqref{rearr1}:
	\begin{equation}
	\label{equationgbyD}
		D\tilde{u} = g.
	\end{equation}
	Let now $w \in \mathcal{D}'(S^{1})$ real-valued be such that $Dw = \tilde{u}$ and $\hat{w}(0) = 0$. Once more, existence of such a distribution $w$ is easily deduced using Fourier coefficients. We would like to emphasise at this point that due to the assumed smoothness of $u, f_j, g_j$, $w$ is smooth as well, as is $F$.\par
	By \eqref{equationgbyD}, we thus notice:
	\begin{equation}
	\label{eqforw}
		D^{2}w = g.
	\end{equation}
	Going over to Fourier coefficients, we see that for every $n \in \mathbb{Z}^\ast$:
	\begin{equation}
		\mathcal{F}\big{(} D^{2}w \big{)}(n) = (1 + i \sign(n))^{2}|n| \hat{w}(n) = 2i \sign(n)|n| \hat{w}(n) = 2in \hat{w}(n) = \hat{g}(n),
	\end{equation}
	or by rearranging:
	\begin{equation}
	\label{fouriercoeffwing}
		\hat{w}(n) = -\frac{i}{2} \frac{\hat{g}(n)}{n}.
	\end{equation}
	Next, we are going to find a suitable distribution $K$ with coefficients $\hat{K}(n) = -\frac{i}{2n}$, in order to express $w$ as a convolution of $g$ with $K$. To this end, let us consider the function $k: [-\pi, \pi] \rightarrow \mathbb{R}$ defined by:
	\begin{equation}
	\label{deffunck}
		k(x)=
\begin{cases}
x+ \pi, \quad \text{if } x < 0 \\
x - \pi, \quad \text{if } x > 0
\end{cases}
	\end{equation}	
	  By slight abuse of notation, let us identify $k$ with its $2\pi$-periodic extension, therefore $k: S^{1} \rightarrow \mathbb{R}$. We calculate the Fourier coefficients of $k$: If $n = 0$, it is obvious due to anti-symmetry that $\hat{k}(0) = 0$. Otherwise, we have $n \neq 0$ and so by using integration by parts:
	\begin{align}
		\hat{k}(n)	&= \frac{1}{2\pi} \int_{-\pi}^{\pi} k(x) e^{-inx}dx \notag \\
				&= \frac{1}{2\pi} \Big{(} \int_{-\pi}^{0} (x + \pi) e^{-inx}dx + \int_{0}^{\pi} (x - \pi) e^{-inx}dx \Big{)} \notag \\
				&= \frac{1}{2\pi} \int_{0}^{\pi} (\pi - x)e^{inx} - (\pi - x)e^{-inx}dx \notag \\
				&= \frac{i}{\pi} \int_{0}^{\pi} (\pi - x) \sin(nx)dx \notag \\
				&= \frac{i}{n} - \frac{i}{\pi} \int_{0}^{\pi} \frac{\cos(nx)}{n} dx = \frac{i}{n}.
	\end{align}
	Consequently, observe that $K = -\frac{1}{2}k$ precisely yields the desired distribution. Let us notice that $K$ is therefore bounded and measurable on $S^{1}$, thanks to the explicit formula for $k$.\\

	Using \eqref{fouriercoeffwing} and the convolution formula for Fourier coefficients, namely:
	\begin{equation}
	\label{convform}
		\widehat{f \ast g}(n) = 2\pi \hat{f}(n) \hat{g}(n), \quad \forall n \in \mathbb{Z},
	\end{equation}
	it is clear that:
	\begin{equation}
	\label{wasconv}
		w =\frac{1}{2\pi} K\ast g.
	\end{equation}
	From \eqref{wasconv}, we obtain by using Young's inequality on $S^{1}$:
	\begin{equation}
	\label{wbounded}
		\| w \|_{L^{\infty}} \leq \frac{1}{2\pi} \| K \|_{L^{\infty}} \| g \|_{L^{1}} = \frac{1}{4} \| g \|_{L^{1}}.
	\end{equation}
	To conclude the first part of the proof, let us observe the following\footnote{This will actually be the first and only point in the proof where we use the fact that $u$ is real-valued in a meaningful way. See Remark \ref{complex} for an extension to complex-valued distributions.}:
	\begin{align}
	\label{simplifynormexpression}
		\int_{S^{1}} (u - F)^{2}dx	&= \int_{S^{1}} Dw (u - F)dx \notag \\
							&\simeq \sum_{n \in \mathbb{Z}} \widehat{Dw}(n) \widehat{u - F}(-n) \notag \\
							&= \sum_{n \in \mathbb{Z}} |n|^{\frac{1}{2}}(1+i \sign(n)) \widehat{w}(n) \widehat{u - F}(-n) \notag \\
							&= \sum_{n \in \mathbb{Z}} \widehat{w}(n) \cdot |n|^{\frac{1}{2}}(1-i \sign(-n)) \widehat{u - F}(-n) \notag \\
							&\simeq \int_{S^{1}} w \overline{D}(u - F)dx \notag \\
							&= \int_{S^{1}} w \overline{D}u dx - \int_{S^{1}} w \overline{D}F dx \notag \\
							&= \int_{S^{1}} w \tilde{g} dx + \int_{S^{1}} w \tilde{f} dx - \int_{S^{1}} w \overline{D}F dx 
	\end{align}
	where we used the Fourier representation of the distribution $u - F$ to justify the second equation. Observe that this enables us to estimate:
	\begin{equation}
	\label{est1}
		\Big{|} \int_{S^{1}} w \tilde{g} dx \Big{|} \leq \| w \|_{L^{\infty}} \| \tilde{g} \|_{L^{1}} \leq \frac{1}{4} \big{(} \| g_{1} \|_{L^{1}} + \| g_{2} \|_{L^{1}} \big{)}^{2},
	\end{equation}
	and:
	\begin{equation}
	\label{est2}
		\Big{|} \int_{S^{1}} w \overline{D}F dx \Big{|} \leq \| w \|_{\dot{H}^{\frac{1}{2}}} \| \overline{D}F\|_{\dot{H}^{-\frac{1}{2}}} \leq C \| u - F \|_{L^{2}} \| f \|_{\dot{H}^{-\frac{1}{2}}}.
	\end{equation}
	The remaining summand may be estimated completely analogous to \eqref{est2}. Notice that we used the explicit definition of the norms of Sobolev spaces with negative exponents and the Fourier multipliers to obtain \eqref{est2}, see \eqref{fouriercoeffest} for the main ideas. Using \eqref{est1} and \eqref{est2} yields:
	\begin{align}
		\| u - F \|_{L^{2}}^{2} 	&\leq \frac{1}{4} \big{(} \| g_{1} \|_{L^{1}} + \| g_{2} \|_{L^{1}} \big{)}^{2} + 2C \| u - F \|_{L^{2}} \| f \|_{\dot{H}^{-\frac{1}{2}}} \notag \\
						&\leq \frac{1}{4} \big{(} \| g_{1} \|_{L^{1}} + \| g_{2} \|_{L^{1}} \big{)}^{2} + \frac{1}{2} \| u - F \|_{L^{2}}^{2} + \frac{4C^{2}}{2} \| f \|_{\dot{H}^{-\frac{1}{2}}}^{2},
	\end{align}
	using the arithmetic-geometric mean inequality. Note that the factor $2C$ is due to estimate \eqref{est2} also applying to the integral of $w\tilde{f}$. By absorbing the $L^{2}$-norm of $u - F$, we arrive at:
	\begin{align}
		\| u - F \|_{L^{2}}^{2}	&\leq \frac{1}{2} \big{(} \| g_{1} \|_{L^{1}} + \| g_{2} \|_{L^{1}} \big{)}^{2} + 4C^{2} \| f \|_{\dot{H}^{-\frac{1}{2}}}^{2} \notag \\
						&\leq \max \{ \frac{1}{2}, 4C^{2} \} \big{(} \| g_{1} \|_{L^{1}} + \| g_{2} \|_{L^{1}} + \| f_{1} \|_{\dot{H}^{-\frac{1}{2}}} + \| f_{2} \|_{\dot{H}^{-\frac{1}{2}}} \big{)}^{2}.
	\end{align}
	Consequently, by estimating the $L^{2}$-norm of $F$ by the $\dot{H}^{-\frac{1}{2}}$-norm of $f_{1}, f_{2}$ using \eqref{fouriercoeffest}, we immediately conclude:
	\begin{equation}
	\label{l2boundu}
		\| u \|_{L^{2}} \leq \tilde{C} \big{(} \| g_{1} \|_{L^{1}} + \| g_{2} \|_{L^{1}} + \| f_{1} \|_{\dot{H}^{-\frac{1}{2}}} + \| f_{2} \|_{\dot{H}^{-\frac{1}{2}}} \big{)}.
	\end{equation}
	The constant $\tilde{C} > 0$ appearing in the estimate is independent of $u, f_{j}, g_{j}$.\\
	
	Now, for a general distribution $u \in \mathcal{D}'(S^{1})$ with $\hat{u}(0) = 0$, let us observe that if we convolute $u$ with a smooth function $\varphi$, the resulting distribution $\varphi \ast u$ will be a smooth function as well (in the sense of regular distributions). By a direct computation, \eqref{condition01} and \eqref{condition02} will continue to hold true if we replace $u, f_{j}, g_{j}$ by their corresponding convolutions with $\varphi$. This is an immediate consequence of the fact that the operators $(-\Delta)^{\frac{1}{4}}, \mathcal{R}$ are Fourier multipliers as well as the linearity of convolutions. Choosing $\varphi$ to be supported on arbitrarily small neighbourhoods of the neutral element in $S^{1}$ (i.e. an approximation of the identity $\varphi_{\varepsilon}$) ensures that the convolutions of $\varphi$ with $f_{j}, g_{j}$ converge in the respective norms as we collapse the support of $\varphi$ (i.e. let the parameter $\varepsilon$ in $\varphi_{\varepsilon}$ tend to $0$) and the approximations of $u$ converge in the distributional sense. As a result, we obtain uniform bounds in the respective spaces. This results in an uniform $L^{2}$-bound for $u$ convoluted with $\varphi_{\varepsilon}$ independent of $\varepsilon$, which can be seen to imply $u \in L^{2}_{\ast}(S^{1})$ by using a weak-$L^2$-convergent subsequence. The estimate follows by the lower semi-continuity of the norm. This concludes our proof.~~\hfill $\Box$\\
	
	\begin{Rem}\label{complex}
	Before we enter the discussion of applications and later a generalisation of Theorem \ref{nonlocalBB1D}, let us quickly discuss the assumption that $u$ is real-valued. In fact, this is merely used at a single point in the proof, namely in \eqref{simplifynormexpression}. However, if we proceed similar to the proof of the generalised result in the next section, i.e. we use:
	\begin{equation}
		\int_{S^1} | u - F|^2 dx = \int_{S^1} (u - F) \cdot \overline{u - F} dx = \int_{S^1} Dw \cdot \overline{u-F} dx \sim \sum_{n \in \mathbb{Z}} \widehat{Dw}(n) \cdot \overline{\widehat{u - F}(n)},
	\end{equation}
	we can easily avoid the use of properties of real-valued distributions. The remainder of the proof follows then completely analogous, i.e. we can remove the assumption of $u$ being real-valued effortlessly. Indeed, this slight generalisation will be key to our applications to Bergman spaces below.
	\end{Rem}
	\section{ Fractional Bourgain-Brezis inequality in the Bergman space ${\mathcal{A}}^2(\D)$}	\label{equiva}

We start with the 
{\bf Proof of Theorem \ref{BBB}.}\\

\noindent
  Let us consider an analytic function $f\colon \D\to \C$  such that   $ \limsup_{r\to 1^-}\|f(re^{i\theta})\|_{L^1+{H}^{-1/2} (S^1)}<+\infty.$  \par
  \noindent
  \textbf{1.} Now let us write $f(z)=\sum_{n\ge 0} f_n z^n$ and
$u(e^{i\theta})=\sum_{n\ge 1}\frac{ f_n}{\sqrt{n}} e^{in\theta}$. We first observe that $f - f(0) = \sum_{n \ge 1} f_n z^n$ and if $f(e^{i \theta}) = g(e^{i\theta}) + h(e^{i\theta})$ with $g \in L^1(S^1)$ and $h \in H^{-1/2}(S^1)$, then:
$$f(e^{i\theta}) - f(0) = g - \dashint_{S^1} g + h - \widehat{h}(0).$$
Note that $h - \widehat{h}(0) \in \dot{H}^{-1/2}(S^1)$ with the norm being controlled by $\| h \|_{H^{-1/2}(S^1)}$. We observe that, using the explicit definitions of the norm:
$$\Big{\|} g - \dashint_{S^1} g \Big{\|}_{L^1(S^1)} + \| h - \widehat{h}(0) \|_{\dot{H}^{-1/2}(S^1)} \lesssim \|g \|_{L^1(S^1)} + \| h \|_{H^{-1/2}(S^1)}$$
Therefore, we may conclude by taking the infimum over all such $g,h$:
\begin{equation}
\label{homogenisationineq}
	\| f - f(0) \|_{L^1 + \dot{H}^{-1/2}(S^1)} \leq C \| f \|_{L^1 + H^{-1/2}(S^1)}.
\end{equation}
Assume therefore first that
$$f(e^{i\theta}) - f(0) =\sum_{n\ge 1}{ f_n} e^{in\theta}\in L^1+\dot{H}^{-1/2}(S^1).$$ 
In this case,  we get 
$(-\Delta)^{1/4}u = f - f(0) \in  L^1+\dot{H}^{-1/2}(S^1) $. Additionally, we observe that since $u$ contains only positive frequencies, we trivially have ${\mathcal{R}}(-\Delta)^{1/4}u\in  L^1+\dot{H}^{-1/2}(S^1)$ as well with
$$\|{\mathcal{R}}(-\Delta)^{1/4}u\|_{L^1+\dot{H}^{-1/2}(S^1) }=\|(-\Delta)^{1/4}u\|_{L^1+\dot{H}^{-1/2}(S^1)}.$$
From the inequality \eqref{fractineq}, observing that $\dashint_{S^1} u   =0$, we deduce that
$$\|u\|_{L^2(S^1)}\le C\|(-\Delta)^{1/4}u\|_ {(L^1+\dot{H}^{-1/2})(S^1)} = C\|f - f(0) \|_ {L^1+\dot{H}^{-1/2}(S^1)} \leq C^\prime \| f \|_ {L^1+{H}^{-1/2}(S^1)},$$
where we used \eqref{homogenisationineq}. Hence $\sum_{n> 0}\frac{{ f_n}}{n} e^{in\theta}\in H^{1/2}(S^1)$ and $g(z)=\sum_{n>0} \frac{f_n}{n} z^n\in H^1(\D)$. We have
$g'(z)=\sum_{n\ge 0} {f_{n+1}} z^{n}\in L^2(\D)$ and 
\begin{eqnarray}\label{bergman4}
\|f(z)-f(0)\|_{L^2(\D)}&=&\|zg'(z)\|_{L^2(\D)}\nonumber\\
&\le& C \|g\|_{H^{1/2}(S^1)}=C\|u\|_{L^2(S^1)}\nonumber\\
&\le& C\|f\|_{L^1+{H}^{-1/2}(S^1)}.\end{eqnarray}
The desired estimate follows by the triangle inequality, if we can show:
$$| f(0) | \leq C \| f \|_{L^1+{H}^{-1/2}(S^1)}$$
To achieve this, let us decompose $f = g + h$ with $g \in L^1(S^1)$ as well as $h \in H^{-1/2}(S^1)$. Then we denote as usual the Fourier coefficients of $g,h$ by $g_n, h_n$ for all $n \in \Z$ and define:
$$G(z) := \sum_{n \geq 0} g_n z^n + \sum_{n < 0} g_n \overline{z}^{|n|}, \quad H(z) := \sum_{n \geq 0} h_n z^n + \sum_{n < 0} h_n \overline{z}^{|n|}.$$
By the summability properties, these define harmonic functions on $\D$ having boundary values $g, h$ respectively. By comparison of the coefficients, we also observe:
$$f(z) = G(z) + H(z),$$
in particular for $z=0$. Moreover, by the mean value property of harmonic functions over the boundary of the disc, we can deduce:
$$| G(0) | \lesssim \| g \|_{L^1(S^1)}.$$
Using the mean value property over the entire disc, we similarily see by H\"older's inequality:
$$| H(0) | \lesssim \| H \|_{L^{1}(\D)} \lesssim \| H \|_{L^{2}(\D)}.$$
It is easy to verify by a direct computation analogous to the same characterisation of the norm in $\mathcal{A}^2(\D)$ that:
$$\| H \|_{L^{2}(\D)}^2 \sim \sum_{n \in \Z} \frac{| h_n|^2}{| n | + 1} \leq \| h \|_{H^{-{1/2}}(S^1)}.$$
In conclusion, we have:
$$| f(0)| \leq C \left( \| g \|_{L^1(S^1)} + \| h \|_{H^{-{1/2}}(S^1)} \right).$$
By taking the infimum over $g,h$ such that $f=g+h$ we get
\begin{equation}\label{bergman4bis}
| f(0)| \leq C \| f \|_{L^1 + H^{-{1/2}}(S^1)}\end{equation}
By combining \eqref{bergman4} and \eqref{bergman4bis}, we obtain the desired estimate:
 \begin{equation}
\label{bergman4tris}
	\|f\|_{L^2(\D)} \leq C\|f\|_{L^1+{H}^{-1/2}(S^1)}.
\end{equation}

\noindent
\textbf{2.} In the general case when   $ \limsup_{r\to 1^-}\|f(re^{i\theta})\|_{L^1+{H}^{-1/2} (S^1)}<+\infty$, we consider
  for every $0<r<1$ the function $f_r(z)=f(rz)\in C^{\infty}(\bar{B}(0,1)).$ We can apply
\eqref{bergman4} to $f_r$ and obtain that
\begin{equation}\label{bergman5}
\|f_r\|_{L^2(\D)}\le C\|f_{r}\|_{ L^1+ {H}^{-1/2}(S^1)}.\end{equation}
Since by assumption $\limsup_{r\to 1^-}\|f(re^{i\theta})\|_{ L^1+ {H}^{-1/2}(S^1)}<+\infty,$
we deduce that
\begin{equation}\label{bergman6}
\sup_{0<r<1}\|f_r\|_{L^2(\D)}<+\infty.
\end{equation}
The inequality implies that actually $f\in L^2(\D)$ as well as\footnote{
Let $f(z)=\sum_{n\ge 0} f_n z^n$
 We observe that 
 \begin{eqnarray}
  \|f_r\|^2_{L^2(\D)}&=&\int_0^1\int_0^{2\pi}|f(\rho re^{i\theta})|^2 \rho d\theta d\rho\nonumber\\
  &=&2\pi\int_0^1\sum_{n=0}^{\infty}|a_n|^2r^{2n}\rho^{2n+1} d\rho=
  2\pi \sum_{n=0}^{\infty}\frac{|a_n|^2r^{2n}}{2n+2}.
  \end{eqnarray}
  and similarily
  \begin{equation}
  \|f\|^2_{L^2(\D)}=2\pi \sum_{n=0}^{\infty}\frac{|a_n|^2}{2n+2}.
  \end{equation}
From \eqref{bergman6} and extracting a weakly convergent subsequence which by convergence of the Fourier coefficients must have limit $f$, it follows that
$\|f\|^2_{L^2(\D)}<+\infty$ and Abel's Theorem on power series yields that
$$
\lim_{r\to 1^-}\|f_r\|^2_{L^2(\D)}=\|f\|^2_{L^2(\D)}.$$
 }
and 
\begin{equation}\label{bergman2}
\|f\|_{L^2(\D)}\le C\|f\|_{ L^1+ {H}^{-1/2}(S^1)}.\end{equation}

\noindent
  {\bf Conversely}, let $f\colon\D\to\C$ be in ${\mathcal{A}}^2(\D^2)$. We write $f(z)=\sum_{n=0}^{\infty} a_nz^n$. We prove the following:\par
  \noindent
{\bf Claim:} $\limsup_{r\to 1^-}\|f(re^{i\theta})\|_{ L^1+{H}^{-1/2}(S^1)}<+\infty$.
\par
\noindent
{\bf Proof of the claim.} We show that $\limsup_{r\to 1^-}\|f(re^{i\theta})\|_{H^{-1/2}(S^1)}<+\infty$. For every
$0<r<1$, we set $f_r(z)=f(rz)\in C^{\infty}(\bar{B}(0,1)).$ Since $f\in L^2(\D)$, we have 
\begin{equation}\label{frf}
\limsup_{r\to 1^-}\|f_r\|_{L^2(\D)}=\|f\|_{L^2(\D)}.
\end{equation}
Moreover
\begin{equation}
\|f_r\|^2_{H^{-1/2}(S^1)}=\sum_{n\ge 0} \frac{|f_n|^2}{1+n} r^{2n}
\end{equation}
and
\begin{eqnarray}\label{h12}
\int_{\D}|f_{r}|^2&=&2\pi\int_0^1\sum_{n\ge0} |f_n|^2 r^{2n}s^{2n} s ds\simeq \sum_{n\ge 0} \frac{|f_n|^2}{2n+2} r^{2n}\nonumber\\
&\simeq& \frac{1}{2} \sum_{n\ge 0} \frac{|f_n|^2}{n+1} r^{2n}\simeq \|f_r\|_{H^{-1/2}(S^1)}.
\end{eqnarray}
By combining \eqref{frf} and  \eqref{h12}, we get that
\begin{equation}
\limsup_{r\to 1^-} \|f_r\|_{ H^{-1/2}(S^1)}\lesssim \|f\|_{L^2(\D)}<+\infty.
\end{equation}
 
  \noindent
We conclude the proof.~~\hfill $\Box$\par

\medskip

Next we show that    Theorem  \ref{BBB} is actually equivalent to Theorem \ref{nonlocalBB1D}.
\begin{prop}\label{equiv}
Theorem \ref{BBB} implies Theorem \ref{nonlocalBB1D}. Therefore, they are equivalent.
 \end{prop}
 \noindent
 {\bf Proof.} We have already seen in the proof of Theorem  \ref{BBB} that Theorem \ref{nonlocalBB1D} implies  the fact  that a holomorphic function with the property that
$\limsup_{r\to 1^-}\|f(re^{i\theta})\|_{ L^1+ {H}^{-1/2}(S^1)}<+\infty$ is in $L^2(\D)$, namely it belongs to the Bergman space ${\mathcal{A}}^2(\D).$

 Conversely, let us consider $u\in C^{\infty}(S^1)$ such that $(-\Delta)^{1/4}u,{\mathcal{R}}(-\Delta)^{1/4}u\in L^1+\dot {H}^{-1/2}(S^1)$. We assume that $\int_{0}^{2\pi} u(e^{i\theta}) d\theta=0$.  We decompose $u=u^++u^-$, where 
  $$u^+=\sum_{n>0}u_n e^{in\theta},~~~u^-=\sum_{n<0}u_n e^{in\theta}.$$
Let us first consider $u^+$. By assumption we have $\sum_{n\ge 1}n^{1/2}u_n e^{in\theta} = 1/2 ((-\Delta)^{1/4}u -i {\mathcal{R}}(-\Delta)^{1/4}u)\in  L^1+\dot {H}^{-1/2}(S^1).$
Let $f(z)=\sum_{n\ge 1}n^{1/2} u_n z^n$ be the harmonic extension of $v=(-\Delta)^{1/4} u^+$ in $\D.$
From Theorem \ref{BBB}, it follows that $f^+=\sum_{n> 0}n^{1/2} u_n z^n\in L^2(\D)$ and
$$\|f^+\|_{L^2(\D)}\le C\|f^+\|_{  L^1+ {H}^{-1/2}(S^1)}\le \|(-\Delta)^{1/4} u^+\|_{  L^1+\dot {H}^{-1/2}(S^1)}.$$
Switching to the homogeneous Sobolev space is possible, as we have $\dot{H}^{-1/2}(S^1) \subset H^{-1/2}(S^1)$ continuously embedded. Since $f^+(z)=\sum_{n>0}n^{1/2} u_n z^n\in L^2(\D)$, it follows that $\sum_{n>0} \frac{u_n}{n^{1/2} }z^n\in   H^1(\D)$ and therefore
 $\sum_{n>0} \frac{u_n}{n^{1/2} }e^{in\theta}\in   \dot{H}^{1/2}(S^1)$.
Hence $u^+\in L^2(S^1)$ with
\begin{eqnarray}
\|u^{+}\|_{L^2(S^1)}&\le& C  \|(-\Delta)^{1/4} u^+\|_{  L^1+\dot {H}^{-1/2}(S^1)}\nonumber \\
&\le& C\left( \|(-\Delta)^{1/4}u\|_{  L^1+\dot {H}^{-1/2}(S^1)}+\|{\mathcal{R}}(-\Delta)^{1/4}u\|_{  L^1+\dot {H}^{-1/2}(S^1)} \right).
\end{eqnarray}
The same arguments hold for $u^-.$ We conclude the proof. ~\hfill $\Box$

	\section{The Bourgain-Brezis Inequality on the Torus ${\mathbf{T^n, n \geq 2}}$}\label{torus}
	
	Next, we would like to generalise the result from Theorem \ref{nonlocalBB1D} to domains of dimension $n \geq 2$. To achieve this while retaining the general structure of the proof, we first have to determine the right set of conditions and the appropriate domain. Observe that it is clear, due to the proof for $S^1$ heavily relying on Fourier series, that the natural domain for such a generalisation is the torus $T^n$. In investigating generalisations of the proof, we have to focus on two aspects: Clifford algebras and boundedness of the kernel. The first subsection introduces complex Clifford algebras as a useful tool in our proof and shows how to generalize the argument presented in section \ref{secproof}. The results and properties of Clifford algebras are due to \cite{gilbert} and \cite{hamilton} and are briefly discussed in the preliminary section \ref{prel}. The second subsection fills in the gaps in the proof by showing that the kernel used is actually bounded, following an argument presented in \cite[p.405-406]{bourgain1}. In the case $n=1$, we have seen that $k$ has an explicit description as a sawtooth function. In higher dimensions, unfortunately, we are not aware of an explicit formula for the kernel. However, due to some estimates on alternating sums, we can remedy this lack of explicit representation and derive the crucial properties abstractly.
	
	 \begin{thm}
	\label{nonlocalBBnD}
		Let $u \in \mathcal{D}'(T^{n})$ be complex-valued and such that:
		$$(-\Delta)^{\frac{n}{4}}u,\mathcal{R}_j(-\Delta)^{\frac{n}{4}}u\in  (L^1 + \dot{H}^{-\frac{n}{2}})(T^n), \quad \forall j \in \{ 1, \ldots n \}.$$
		  Then we have 
		 			$u - \dashint_{T^{n}} u dx \in L^{2}_{\ast}(T^{n})$ with
		\begin{equation}
		\label{estimateforuinthm}
			\Big{\|} u - \dashint_{T^{n}} u dx \Big{\|}_{L^{2}} \leq C\left(\|(-\Delta)^{n/4} u\|_{L^1 + \dot{H}^{-n/2}(T^n)}+\sum_{j=1}^n\|{\mathcal{R}}_j(-\Delta)^{n/4} u\|_{L^1 + \dot{H}^{-n/2}(T^n)}\right),
\end{equation}
		for some $C > 0$ independent of  $u$.
	\end{thm}
\noindent	
{\bf Proof of Theorem \ref{nonlocalBBnD}.}\par
	
	Let us first note that, if we take $T^n$ with $n \geq 2$, there are $n$ different Riesz transforms, one for each basis direction. This suggests that the right conditions should involve some restriction on each of the Riesz transforms. In addition, considering the symbol of $D^2$, we see that we rely some cancellation property stemming from the complex nature of $i$. \footnote{This refers to the property $i^2 = -1$ which was key to reduce the multiplier of $D^2$ to a simpler form.}
	Therefore, a natural way to obtain a generalisation would involve Clifford algebras to include sufficiently many anticommuting complex units.  	
	\par
			Firstly, it is immediate that the same simplifications as in the case $n=1$ apply here. So we may assume $\hat{u}(0) = 0$. Throughout most of this proof, the coefficient $m=0$ will be implicitly omitted, as it will be vanishing for all functions/distributions considered. Moreover, the reduction to smooth functions applies equally well in this case. Therefore, we may assume without loss of generality that  $u, f_j, g_j$ are all smooth.\\
		
		The heart of the argument lies in the correct definition of $D$ and $\overline{D}$ on $T^n$. As mentioned in the introduction of the current section, Clifford algebras and their set of complex units actually provide the desired framework. Let $\mathbb{C}_n$ denote the universal complex Clifford algebra associated with the quadratic space $(\mathbb{C}^{n}, Q)$,
		where:
		\begin{equation}
			Q(z_1, \ldots, z_n) := -\sum_{j=1}^{n} z_j^2, \quad \forall (z_1, \ldots, z_n) \in \mathbb{C}^{n}. 
		\end{equation}
		We emphasise that the particular choice of $Q$ is at odds with usual conventions for complex Clifford algebras, but using our quadratic form, we obtain the appropriate basis commutation relations while remaining isomorphic to the usual convention. One could reduce to the usual defining quadratic form by choosing $i \cdot e_j$ instead of the standard basis $e_j$ throughout our proof. In fact, the main reason why we decided to use our convention is to use the Riesz operators in their usual form.\\
		Observe that we then have, for the standard basis denoted by $e_1, \ldots, e_n$:
		\begin{equation}
			e_j e_k + e_k e_j = 2 \delta_{jk}, \quad \forall j, k \in \{1, \ldots, n\},
		\end{equation}
		simply by the definition of Clifford algebras and the quadratic form $Q$. We define now for any $v\in C^{\infty}(T^n, \mathbb{C})$:
		\begin{align}
		\label{diffop1}
		Dv 			&= \Delta^{\frac{n}{4}}( Id + \sum_{j=1}^{n} e_j \mathcal{R}_{j})v \\
		\label{diffop2}
		\overline{D}v 	&= \Delta^{\frac{n}{4}}( Id - \sum_{j=1}^{n} e_j \mathcal{R}_{j})v.
		\end{align}
		We emphasise the similarity with \cite[(5.14)]{gilbert} used in the context of Hardy spaces. The crucial observation for our purposes is the following multiplier property for Fourier series for every $m \in \mathbb{Z}^{n}$:
		\begin{align}
		\mathcal{F}(Dv)(m) 				&= | m |^{\frac{n}{2}} \big{(} 1 + \sum_{j=1}^{n} e_j \cdot i \frac{m_j}{| m |} \big{)} \mathcal{F}(v)(m) \\
		\mathcal{F}(\overline{D}v)(m)	 	&= | m|^{\frac{n}{2}} \big{(} 1 - \sum_{j=1}^{n} e_j \cdot i \frac{m_{j}}{| m |} \big{)} \mathcal{F}(v)(m),
		\end{align}
		where $| m|$ denotes the Euclidean norm on $\mathbb{Z}^{n}$. We highlight that at this point, we know that $Du$ and $\overline{D}u$ are functions in $L^1 + \dot{H}^{-\frac{n}{2}}(T^n, \mathbb{C}_{n})$. Completely analogous to the proof of Theorem \ref{nonlocalBB1D}, we may find $F \in L^{2}$ (due to the invertibility of non-zero vectors $v \in \mathbb{R}^{n}$ in $\mathbb{C}_{n}$\footnote{Observe that for real vectors in $\mathbb{R}^n$, we find $m^2 = |m|^2$. For general vectors in $\mathbb{C}^{n}$, this fails, as can be seen in the counterexample: $$(e_1 + i e_2 )^2 = 0$$}). To be precise, observe that if $DF = f$, $f$ and $g$ are defined to satisfy $Du = f + g$ by splitting the terms $f_j, g_j$ in the natural way, then:
		\begin{equation}
			\forall m \neq 0: \quad | m |^{\frac{n}{2}} \big{(} 1 + \sum_{j=1}^{n} e_j \cdot i \frac{m_j}{| m |} \big{)} \widehat{F}(m) = \widehat{f}(m),
		\end{equation}
		which may be rewritten as:
		\begin{equation}
			\widehat{F}(m) = \frac{1}{2|m|^{\frac{n}{2}}} \Big{(} 1 - \sum_{j=1}^{n} e_j \cdot i \frac{m_j}{|m|} \Big{)} \widehat{f}(m),
		\end{equation}
		by using the multiplication relations and associativity on $\mathbb{C}_{n}$. To conclude that $F \in L^{2}$, it suffices to check summability of the Fourier coefficients:
		\begin{align}
			\sum_{m \in \mathbb{Z}^{n} \setminus \{ 0 \}} \| \widehat{F}(m) \|^{2}	&= \sum_{m \neq 0} \Big{\|} \frac{1}{2|m|^{\frac{n}{2}}} \big{(} 1 - \sum_{j=1}^{n} e_j \cdot i \frac{m_j}{|m|} \big{)} \widehat{f}(m) \Big{\|}^2 \notag \\
																&\lesssim \sum_{m \neq 0} \frac{1}{|m|^{n}} \| \widehat{f}(m) \|^2 \notag \\
																&\lesssim \| f \|_{\dot{H}^{-\frac{n}{2}}}^2 < +\infty.
		\end{align}
		We mention here that the characterisations for regularity and integrability carry over without problem, even if we use Clifford algebra-valued functions by verifying componentwise regularity.\\
		
		Consequently, as in the case $n=1$, we may define $\tilde{u} = u - F$ and observe that $D\tilde{u} =: g \in L^1$. Solving $Dw = \tilde{u}$ in the sense of distributions leaves us with $D^{2}w = g$.\\
		
		The key point behind the second proof of Theorem \ref{nonlocalBB1D}  lies in the fact, that $D^2$ has an inverse given by the convolution with a bounded function. By a direct computation, we arrive at the following expression for the multiplier associated with $D^2$:
		\begin{equation}
		\label{intermediate1}
			\mathcal{F}(D^2w)(m) = | m |^{n} \big{(} 1 + \sum_{j=1}^{n} e_j \cdot i \frac{m_j}{| m |} \big{)}^2 \mathcal{F}(w)(m) = 2i \cdot | m |^{n} \big{(} \sum_{j=1}^{n} e_j \frac{m_j}{| m |} \big{)} \mathcal{F}(w)(m),
		\end{equation}
		for every $m \in \mathbb{Z}^{n}$. Observe that we used the fact that the complex unit $i$ of $\mathbb{C}$ commutes with all $e_j$ (as the Clifford algebra is a complex algebra)
		and that:
		\begin{equation}
			(i \cdot e_j)^2 = i^2 \cdot e_j^2 = i^2 = -1.
		\end{equation}
		Let us identify $m = \sum m_j e_j$, i.e. we consider the vector $m \in \mathbb{Z}^{n} \subset \mathbb{C}^{n}$ as an element in $\mathbb{C}_n$. Therefore, \eqref{intermediate1} becomes:
		\begin{equation}
			\mathcal{F}(D^2w)(m) = 2i | m |^{n-1} \cdot m \mathcal{F}(w)(m), \quad \forall m \in \mathbb{Z}^{n}.
		\end{equation}
		As stated before, all vectors $\mathbb{R}^{n} \subset \mathbb{C}_n$
		 are invertible due to:
		\begin{equation}
			z^2 = -Q(z), \quad \forall z \in \mathbb{C}^{n}.
		\end{equation}
		So, for the real vector $m$, we have due to $m \cdot m = - Q(m)$:
		\begin{equation}
			m^{-1} = \frac{m}{| m |^2}, \quad \forall 0 \neq m \in \mathbb{Z}^{n}.
		\end{equation}
		This means that $D^2 w = g$ can be restated as:
		\begin{equation}
			\mathcal{F}(w)(m) = \frac{1}{2i} \cdot \frac{m}{| m |^{n+1}} \mathcal{F}(g)(m),
		\end{equation} 
		for every $0 \neq m \in \mathbb{Z}^{n}$.\\
		
		For now, let us assume that a function $K$ on the torus exists, such that:
		\begin{equation}
			\widehat{K}(m) = \frac{1}{2i} \cdot \frac{m}{| m |^{n+1}}, \quad \forall m \in \mathbb{Z}^{n} \setminus \{0\}.
		\end{equation}
		In this case, we may check using Fourier coefficients that (keeping in mind that the order of factors in the convolution matters for products in Clifford algebras):
		\begin{equation}
			w = \frac{1}{(2\pi)^n} K \ast g 
		\end{equation}
		Thus, we have the following inequality:
		\begin{equation}
			\| w \|_{L^{\infty}} \lesssim \| K \|_{L^{\infty}} \| g \|_{L^{1}}.
		\end{equation}
		This is an immediate consequence of the definition, Minkowski's inequality and continuity of the Clifford multiplication in the Clifford algebra norm.\\
		Moreover, we may deduce:
		\begin{align}
			\| u - F \|_{L^2}^2	&= \int_{T^{n}} P_{0}\big{(} \overline{(u-F)} \cdot (u-F) \big{)}dx \notag \\
							&= P_{0} \Big{(} \int_{T^{n}} \overline{(u-F)} \cdot (u-F) dx \Big{)} \notag \\
							&\leq \Big{\|} \int_{T^{n}} \overline{(u-F)} \cdot (u-F) dx \Big{\|} \notag \\
							&= \Big{\|} \int_{T^{n}} \overline{Dw} \cdot (u-F) dx \Big{\|} \notag \\
							&= \Big{\|}\int_{T^{n}} \Big{(} \sum_{m} \overline{\widehat{Dw}(m)} e^{-i\langle m, x \rangle} \Big{)} \cdot \Big{(}\sum_{\tilde{m}} \widehat{u-F}(m) e^{i \langle \tilde{m}, x \rangle} \Big{)} dx \Big{\|} \notag \\
							&\simeq \Big{\|} \sum_{m} \overline{\widehat{Dw}(m)} \cdot \widehat{u-F}(m) \Big{\|} \notag \\
							&= \Big{\|} \sum_{m} \overline{| m |^{\frac{n}{2}} \big{(} 1 + \sum_{j=1}^{n} e_j \cdot i \frac{m_j}{| m |} \big{)}\widehat{w}(m)} \cdot \widehat{u-F}(m) \| \notag \\
							&= \Big{\|} \sum_{m} \overline{\widehat{w}(m)} \cdot \overline{| m |^{\frac{n}{2}} \big{(} 1 + \sum_{j=1}^{n} e_j \cdot i \frac{m_j}{| m |} \big{)}}\widehat{u-F}(m) \| \notag \\
							&= \Big{\|} \sum_{m} \overline{\widehat{w}(m)} \cdot | m |^{\frac{n}{2}} \big{(} 1 - \sum_{j=1}^{n} e_j \cdot i \frac{m_j}{| m |} \big{)}\widehat{u-F}(m) \| \notag \\
							&= \Big{\|} \sum_{m} \overline{\widehat{w}(m)} \cdot \widehat{{\overline{D}}(u-F)}(m) \| \notag \\
							&\simeq \Big{\|} \int_{T^{n}} \overline{w} \cdot {\overline{D}}(u-F) dx \Big{\|}.
		\end{align}
		Observe that in the first inequality, we used that the norm squared of $u - F$ actually appears as the coefficient associated with $1$ in the product $\overline{u-F} \cdot (u-F)$. In addition, the conjugation in the ninth line can easily deduced from our definition in the preliminary section of the paper, see \eqref{cliffconjdef}. The remainder of the argument then follows completely analogous to the $1D$-proof, up to the obvious modifications. Again, simple considerations show that we even have the following inequality:
		\begin{equation}
		\label{ineqonsum}
			\Big{\|} u - \dashint_{T^{n}} u dx \Big{\|}_{L^{2}} \lesssim \sum_{j=0}^{n} \big{\|} (-\Delta)^{\frac{n}{2}} \mathcal{R}_{j} u \big{\|}_{\dot{H}^{-\frac{n}{2}} + L^{1}},
		\end{equation}
		where $\mathcal{R}_{0} = Id$.\\
		
		
		To complete the proof in the same way as for Theorem \ref{nonlocalBB1D}, we still need to find a bounded kernel $K$ satisfying:
		\begin{equation}
			\widehat{K}(m) = \frac{1}{2i} \cdot \frac{m}{| m |^{n+1}}, \quad \forall 0 \neq m \in \mathbb{Z}^{n}.
		\end{equation}
		This is the purpose of the next subsection, so we may conclude the proof of Theorem \ref{nonlocalBBnD} at this point.~~\hfill $\Box$
	
	\subsection{Boundedness of the Kernel}
	
	Lastly, let us find an appropriate kernel. We first notice that due to linearity, symmetry and the splitting into different directions, it is enough to find a bounded function $k$, such that:
	\begin{equation}
		\widehat{k}(m) = \frac{m_1}{| m |^{n+1}}, \quad \forall 0 \neq m \in \mathbb{Z}^{n}.
	\end{equation}
	Consequently, we want to study the boundedness of the following conditionally convergent series:
	\begin{equation}
	\label{formulakernel1}
		k(x) = \sum_{m \in \mathbb{Z}^{n} \setminus \{ 0\}} \frac{m_1}{| m |^{n+1}} e^{i \langle m, x \rangle}.
	\end{equation}
	Let us fix some notation. We usually identify $m \in \mathbb{Z}^{n}$ with $m = (m_1, \tilde{m})$, where $\tilde{m} \in \mathbb{Z}^{n-1}$. We will sometimes use the same notation for $x \in \mathbb{R}^{n}$. Moreover, for any $m$, we define $m^\prime = (-m_1, \tilde{m})$. This allows us to immediately see:
	\begin{equation}
		\widehat{k}(m^\prime) = - \widehat{k}(m), \quad \forall m \in \mathbb{Z}^{n} \setminus \{ 0 \}.
	\end{equation}
	This observation enables us to rewrite \eqref{formulakernel1} as follows:
	\begin{equation}
		\label{formulakernel2}
		k(x) = 2i \cdot \sum_{m_1 > 0} \sum_{\tilde{m} \in \mathbb{Z}^{n-1}} \frac{m_1}{| m |^{n+1}} \sin(m_1 x_1) e^{i \langle \tilde{m}, \tilde{x} \rangle}.
	\end{equation}
	The strategy of the proof is based on \cite[p.405-406]{bourgain1}. Thus, the main point is to split the sum into partial sums involving $m_1$ and $| \tilde{m} |$ being comparable to some dyadic $2^{k_1}$ and $2^{\tilde{k}}$ respectively. Then, we distinguish $k_1 \leq \tilde{k}$ and $k_1 \geq \tilde{k}$ to conclude. Thus, we consider the following sum derived from \eqref{formulakernel2}:
	\begin{equation}
		\label{estimatekernel1}
		| k(x) | \leq \sum_{k_1 \geq 0} \sum_{\tilde{k} \geq 0} \Big{|} \sum_{m_1 \sim 2^{k_1}} \sum_{| \tilde{m} | \sim 2^{\tilde{k}}} \frac{m_1}{| m |^{n+1}} \sin(m_1 x_1) e^{i \langle \tilde{m}, \tilde{x} \rangle} \Big{|}.
	\end{equation}
	Let us mention an uniform estimate for fixed $k_1, \tilde{k}$. To achieve this, we distinguish two cases: $k_1 \geq \tilde{k}$ and $k_1 < \tilde{k}$. We shall need the following estimate that can be found in \cite[(4.22)]{bourgain1}:
	\begin{equation}
	\label{sinest}
		\Big{|} \sum_{\ell \in I} \sin(\ell x) \Big{|} \lesssim 4^k |x| \wedge \frac{1}{|x|},
	\end{equation}
	for every $k \in \mathbb{N}$, $x \in S^1$ and subinterval $I \subset [2^{k - 1}, 2^{k}]$. Here, $\wedge$ denotes the minimum of two functions. Let us provide the argument in a more abstract manner: Consider a finite sum of the form:
	\begin{equation}
	\label{abstractseries}
		\sum_{m_1} \sum_{\tilde{m}} a_{m_1} b_{\tilde{m}} c_{m_1, \tilde{m}}.
	\end{equation}
	Observe that the summands in \eqref{estimatekernel1} inside the absolute value clearly have this form. Let us denote by $A_{m_1}$ the partial sum of all $a_l$ up to the $m_1$-th element. In the case of \eqref{estimatekernel1}, this would be a sum of $\sin(l x)$ over an interval with $l$ comparable to $2^{k_1}$, hence we may use the bound \eqref{sinest}. Therefore, we may rewrite \eqref{abstractseries} as:
	\begin{align}
	\label{simplification}
		\sum_{m_1} \sum_{\tilde{m}} a_{m_1} b_{\tilde{m}} c_{m_1, \tilde{m}}	&= \sum_{m_1} \sum_{\tilde{m}} (A_{m_1} - A_{m_1 - 1} ) b_{\tilde{m}} c_{m_1, \tilde{m}} \notag \\
																&= \sum_{m_1} \sum_{\tilde{m}} A_{m_1} b_{\tilde{m}} (c_{m_1, \tilde{m}} - c_{m_1+1, \tilde{m}}),
	\end{align}
	which, in the case of \eqref{estimatekernel1}, can be estimated using the bound on sums of sinus functions in \eqref{sinest}, the boundedness of the $b_{\tilde{m}}$ which are merely $e^{i\langle \tilde{m}, \tilde{x} \rangle}$ and finally the estimate:
	\begin{equation}
	\label{differentiationest}
		\Big{|}  \frac{m_1}{| m |^{n+1}}  -  \frac{m_1 + 1}{( (m_1 + 1)^2 + | \tilde{m} |^{2})^{\frac{n+1}{2}}} \Big{|} \lesssim \frac{1}{| m |^{n+1}}.
	\end{equation}
	We mention the slight imprecision, as in \eqref{simplification}, the extremal partial sums $A_l$ require further attention. However, in the case we are considering, similar techniques can be applied (since we no longer sum over $m_1$) and we omit further details.\\
	
	Therefore, we arrive at the following estimate:
	\begin{equation}
	\label{baseestimate}
		\Big{|} \sum_{m_1 \sim 2^{k_1}} \sum_{| \tilde{m} | \sim 2^{\tilde{k}}} \frac{m_1}{| m |^{n+1}} \sin(m_1 x_1) e^{i \langle \tilde{m}, \tilde{x} \rangle} \Big{|} \lesssim 2^{k_1}\big{(} 2^{k_1} |x_1| \wedge \frac{1}{2^{k_1}|x_1|} \big{)} \Big{\|} \frac{1}{| m |^{n+1}} \Big{\|}_{l^1(m_1 \sim 2^{k_1}, | \tilde{m} | \sim 2^{\tilde{k}})}
	\end{equation}
	If $k_1 \geq \tilde{k}$, we may simplify \eqref{estimatekernel1} using \eqref{baseestimate} as follows:
	\begin{align}
		| k(x)| 	&\lesssim \sum_{k_1 \geq 1} 2^{\tilde{k}(n-1)} 2^{k_1} \frac{1}{2^{k_1 (n+1)}} \cdot 2^{k_1}\big{(} 2^{k_1} |x_1| \wedge \frac{1}{2^{k_1}|x_1|} \big{)} \notag \\
				&\leq \sum_{k_1 \geq 0} 2^{k_1} |x_1| \wedge \frac{1}{2^{k_1}|x_1|} \lesssim C < \infty, 
	\end{align}
	which can be easily bounded by the definition of the minimum.\\
	
	\noindent
	If $\tilde{k} > k_1$, we find:
	\begin{align}
		\sum_{k_1 \geq 0} \sum_{\tilde{k} \geq 0} &\Big{|} \sum_{m_1 \sim 2^{k_1}} \sum_{| \tilde{m} | \sim 2^{\tilde{k}}} \frac{m_1}{| m |^{n+1}} \sin(m_1 x_1) e^{i \langle \tilde{m}, \tilde{x} \rangle} \Big{|} \notag \\	
					&\lesssim \sum_{k_1} \sum_{\tilde{k} > k_1} 4^{k_1}\big{(} 2^{k_1} |x_1| \wedge \frac{1}{2^{k_1}|x_1|} \big{)} \cdot \frac{1}{2^{k_1 (n+1)}} \sum_{| \tilde{m} | \sim 2^{\tilde{k}}} \frac{1}{(1 + \frac{| \tilde{m}|^2}{2^{2 k_1}})^{\frac{n+1}{2}}} \notag \\
					&\lesssim \sum_{k_1 \geq 0} \frac{2^{k_1 (n-1)} 4^{k_1}}{2^{k_1 (n+1)}} \big{(} 2^{k_1} |x_1| \wedge \frac{1}{2^{k_1}|x_1|} \big{)} \notag \\
					&\leq \sum_{k_1 \geq 0} 2^{k_1} |x_1| \wedge \frac{1}{2^{k_1}|x_1|} \leq C < \infty,
	\end{align}
	where we estimated the sum over $\tilde{m}, \tilde{k}$ by a dominating integral. This shows that $k(x)$ is actually bounded and possesses the required Fourier coefficients, hence adding the last ingredient missing in our proof of Theorem \ref{nonlocalBBnD}.
	

\section{Existence Result for a certain Fractional PDE}\label{exist}

Similar to \cite{bourgain1}, the estimates in Theorem \ref{nonlocalBB1D} and \ref{nonlocalBBnD} may be used to derive existence results for a particular differential operator. However, before turning to the PDE itself, let us briefly provide an alternative formulation of our main theorems for a more general class of distributiions:
\begin{thm}
\label{mainres3}
	Let $u \in \mathcal{D}^\prime(T^{n}, \mathbb{C}_n)$ be $\mathbb{C}_n$-valued and assume that:
	\begin{equation}
		Du, \overline{D}u \in \dot{H}^{-\frac{n}{2}} + L^{1}(T^{n}, \mathbb{C}_{n}).
	\end{equation}
	Here, $D$ and $\overline{D}$ are the operators defined in the proof of Theorem \ref{nonlocalBBnD}. Then $u \in L^{2}(T^n, \mathbb{C}_n)$ and we have the following estimate:
	\begin{equation}
		\Big{\|} u - \int_{T^{n}} u dx \Big{\|}_{L^2} \lesssim \| Du \|_{\dot{H}^{-\frac{n}{2}} + L^{1}} + \| \overline{D}u \|_{\dot{H}^{-\frac{n}{2}} + L^{1}}.
	\end{equation}
\end{thm}
This result is an immediate corollary of the proof of Theorem \ref{nonlocalBBnD}, as we always work with $Du$ and $\overline{D}u$ rather than the $\mathcal{R}_{j} (-\Delta)^{\frac{n}{2}}$. The possibility to generalise to Clifford algebra-valued distributions follows directly, as all arguments involved behave well with respect to the Clifford algebra product. One could also rewrite the estimate by separating the identity operator from the Riesz operators.\\

Let us now turn to the existence result. We would like to consider the following problem:
\begin{equation}
\label{decompositionofg}
	g =  (-\Delta)^{\frac{n}{4}}f_0+\sum_{j=1}^{n} (-\Delta)^{\frac{n}{4}} \bar{\mathcal{R}}_{j} f_{j},
\end{equation}
where   $g \in L^{2}_{\ast}(T^{n})=\lf\{u\in L^2(T^n):~~\dashint_{T^n} u=0\rg\}$.\footnote{The conjugate operator $\bar{\mathcal{R}}_j$ appears due to the duality used in the proof. This ensures, that we can apply the result in Theorem \ref{nonlocalBBnD}. It is simpel to see that by suitably exchanging $\mathcal{R}_j$ by $\bar{\mathcal{R}}_j$ throughout the proof of Theorem \ref{nonlocalBBnD}, the same inequality can be obtained for the dual operators and thus yields the same result as in Corollary \ref{cor1} for the usual Riesz operators.} Obviously, the PDE admits solutions $f_{0}, \ldots, f_{n}$ in $\dot{H}^{\frac{n}{2}}(T^n)$. Again, using Sobolev embeddings, it is also clear that there is a-priori no way to deduce that the $f_j$ may be chosen to be bounded or even continuous. We shall remedy this apparent lack of regularity:

\begin{cor}
\label{cor1}
	Let $g \in L^{2}_{\ast}(T^{n})$. Then there exist $f_{0}, \ldots, f_{n} \in \dot{H}^{\frac{n}{2}} \cap C^{0}(T^{n})$, such that \eqref{decompositionofg} holds.
\end{cor}

\noindent{\bf Proof of Corollary \ref{cor1}.}
	The proof is completely analogous to the one in \cite[Proof of Theorem 1]{bourgain1}: Let us define the following operator:
	\begin{equation}
		T: \bigoplus_{j=0}^{n} \dot{H}^{\frac{n}{2}} \cap C^{0}(T^{n}) \rightarrow L^{2}_{\ast}(T^{n}), \quad T(u_0, \ldots, u_n) := (-\Delta)^{\frac{n}{4}}u_0+\sum_{j=1}^{n} (-\Delta)^{\frac{n}{4}} \bar{\mathcal{R}}_{j} u_{j}.	\end{equation}
	 	It is clear that $T$ is a bounded, linear operator. Moreover, we have that its dual operator is given by:
	\begin{equation}
		T^{\ast}: L^{2}_{\ast}(T^{n}) \rightarrow \bigoplus_{j=0}^{n} \dot{H}^{-\frac{n}{2}} + \mathcal{M}(T^{n}), \quad T^{\ast}(v) := \big{(} (-\Delta)^{\frac{n}{4}}   v, \mathcal{R}_{1}(-\Delta)^{\frac{n}{4}}   v,  \ldots, \mathcal{R}_{n}(-\Delta)^{\frac{n}{4}}   v \big{)}.
	\end{equation}
	Here, $\mathcal{M}(T^{n})$ denotes the collection of Radon measures on $T^n$. As in \cite[(4.3)]{bourgain1}, it can be easily seen (using convolutions) that:
	\begin{equation}
		\| \cdot \|_{\dot{H}^{-\frac{n}{2}} + \mathcal{M}} = \| \cdot \|_{\dot{H}^{-\frac{n}{2}} + L^{1}} \quad \text{on } \dot{H}^{-\frac{n}{2}} + L^{1}(T^{n}).
	\end{equation}
	Therefore, we know by \eqref{ineqonsum} that:
	\begin{equation}\label{ineqonsum2}
		\| u \|_{L^{2}} \lesssim \| T^{\ast} u \|_{\bigoplus \dot{H}^{-\frac{n}{2}} + \mathcal{M}(T^{n})}.
	\end{equation}
	This implies that $T$ is surjective (see Theorem 2.20 in \cite{Brez}). The open mapping Theorem yields that
there is $C>0$ such that
$B^{L^2_*}(0,C)\subseteq T(B^{E}(0,1))$, where $E=\bigoplus_{j=0}^{n} \dot{H}^{\frac{n}{2}} \cap C^{0}(T^{n}).$
 Therefore, for every $g\in L^2_*(S^1)$, there are $(f_0,\ldots,f_n)\in E$ such that
$(-\Delta)^{1/4} f_0+\sum_{i=1}^n(-\Delta)^{1/4} \bar{\mathcal{R}}_j f_j=g$ and 
\begin{equation}
		\sum_{j=0}^{n} \| f_{j} \|_{\dot{H}^{\frac{n}{2}} \cap L^{\infty}} \leq C \| g \|_{L^2},
	\end{equation}
	for   some fixed $C > 0$. This concludes the proof.~~\hfill $\Box$\\

Using Corollary \ref{cor1}, we may derive the following simple result:

\begin{cor}
\label{cor2}
	Let $f \in \dot{H}^{\frac{n}{2}}(T^{n})$. Then there exist $f_{0}, \ldots, f_{n} \in \dot{H}^{\frac{n}{2}} \cap C^{0}(T^{n})$ as well as a smooth function $\varphi \in C^{\infty}(T^{n})$, such that:
	\begin{equation}
		f = \varphi + \sum_{j=0}^{n} \mathcal{R}_{j} f_{j}.
	\end{equation}
\end{cor}

\noindent{\bf Proof of Corollary \ref{cor2}.}
	Take $g = (-\Delta)^{\frac{n}{4}}f \in L^{2}_{\ast}(T^{n})$. By Corollary \ref{cor1}, we see that there exist $f_{0}, \ldots, f_{n} \in \dot{H}^{\frac{n}{2}} \cap C^{0}(T^{n})$, such that \eqref{decompositionofg} is satisfied. Therefore, we know:
	\begin{equation}
		(-\Delta)^{\frac{n}{4}} \Big{(} f -  \sum_{j=0}^{n} \mathcal{R}_{j} f_{j} \Big{)} = 0.
	\end{equation}
	But this implies that the difference lies in the kernel of $(-\Delta)^{m}$, where $m$ is the smallest integer larger or equal than $\frac{n}{4}$. Thus the difference is smooth, leading to the desired decomposition.~~\hfill $\Box$

\section{Appendix}\label{appendix}
In this section, we provide for the reader's convenience a proof of the two inequalities \eqref{embHA} and \eqref{embHA2}, 
	since the authors have not found a precise reference in the literature. \\
	
	\par
	\noindent
	{\bf 1.} Assume first that $f(z)=\sum_{n\ge 0} a_nz^n$ is an analytic function such that $\lim_{r\to 1^-}\|f(re^{i\theta})\|_{L^1(S^1)}<+\infty$.
	 Let $h\in L^{1}(S^1)$ be such that $\lim_{r\to 1^-}\|f(re^{i\theta})-h\|_{L^1(S^1)}.$
	We set $g(z)=\sum_{n\ge 0}\frac{ a_n}{n+1}z^{n+1}$. We observe that $g^{\prime}(z)=f(z)$. From our hypothesis, we have 
	 $ \lim_{r\to 1^-}\|g'(re^{i\theta})\|_{L^1(S^1)}<+\infty.$ Observe that this implies that 
	$\lim_{r\to 1^-} (\|\partial_{\theta}g(re^{i\theta})\|_{L^1(S^1)}+\|\partial_{r}g(re^{i\theta})\|_{L^1(S^1)})<+\infty$. 
	Define $g_r(z)=g(rz)$ for $0<r<1$. Since $g$ is harmonic in  $\D$,
	  we have
	\begin{eqnarray}\label{estL1}
	0&=&\int_{\D}(\Delta g_r \bar g_r+g_r\Delta \bar g_r )dx =\int_{\partial\D}(\partial_r g_r\cdot \bar g_r+ \partial_r\bar g_r g_r)d\sigma -2 \int_{\D}|\nabla g_r|^2 dx\nonumber\\
	&=& \int_{\partial\D}(\partial_r g\cdot \bar g+ \partial_r\bar g g)d\sigma -\int_{\D}|g^{\prime}_r|^2 dx.
	\end{eqnarray}
	We first have (observe that $\dashint_{S^1} g_r=0$)
	\begin{equation}\label{estLinfty}
	\|g_r\|_{L^{\infty}{(S^1)}}\lesssim \|\partial_{\theta} g_r\|_{L^1(S^1)}
	\end{equation}
	and
	from \eqref{estL1} it follows that
 	\begin{equation}\label{estL2}
	\|f_r\|_{L^2(\D)}\simeq\|g^{\prime}_r\|_{L^{2}{(\D)}}\lesssim \|g_r\|_{L^{\infty}{(S^1)}}\|\partial_{r} g\|_{L^1(S^1)}\lesssim  \|  g_r'\|^2_{L^1(S^1)}.
	\end{equation}
	We let $r\to 1$ in \eqref{estL2} and get
	\begin{equation}\label{estL3}
	\|f\|_{L^2(\D)} \lesssim \|h\|_{L^1(S^1)}.
	\end{equation}
	
	\noindent
	{\bf 2.} Assume now that $f(z)=\sum_{n\ge 0} a_nz^n$ is an analytic function such that: 
	$$\lim_{r\to 1^-}\|f(re^{i\theta})\|_{H^{-1/2}(S^1)}<+\infty.$$\par
	
	\noindent
	{\bf Claim.} Assume $a_0 = 0$ in the power series above. Then the series $\sum_{n \ge 1} \frac{|a_n|^2}{n}<+\infty$ and 
	$$\sum_{n\ge 1} \frac{|a_n|^2}{n}=\lim_{r\to 1^-}\|f(re^{i\theta})\|^2_{\dot H^{-1/2}(S^1)}.$$
	\noindent
	{\bf Proof of the claim.} We set $A=\lim_{r\to 1^-}\|f(re^{i\theta})\|^2_{H^{-1/2}(S^1)} \simeq \lim_{r\to 1^-}\|f(re^{i\theta})\|^2_{\dot H^{-1/2}(S^1)}.$ We observe that: 
	$$\|f(re^{i\theta})\|^2_{\dot{H}^{-1/2}(S^1)}=\sum_{n >0} \frac{|a_n|^2r^{2n}}{n}.$$
	For every $N>1$, we have
	\begin{eqnarray}
	A\ge \lim_{r\to 1^-} \sum_{n=1}^N  \frac{|a_n|^2r^{2n}}{n} &=&\sum_{n=1}^N  \lim_{r\to 1^-}\frac{|a_n|^2r^{2n}}{n}\nonumber\\
	&=&\sum_{n=1}^N \frac{|a_n|^2}{n}.
	\end{eqnarray}
	By letting $N\to +\infty$, we get $\sum_{n=1}^\infty \frac{|a_n|^2}{n}<+\infty$ and by Abel's theorem on power series, we deduce that the norms converge
	$$
	 \lim_{r\to 1-}\sum_{n> 0} \frac{|a_n|^2r^{2n}}{n}=\sum_{n=1}^\infty \frac{|a_n|^2}{n}.$$
	Therefore, $f(e^{i\theta})\in \dot{ H}^{-1/2}(S^1)$ 
	and $\lim_{r\to 1^-}\|f(re^{i\theta})-f(e^{i\theta})\|_{\dot H^{-1/2}(S^1)}=0$, by observing that the convergence holds weakly and the norms converge, which is an equivalent characterisation for convergence with respect to the norm in Hilbert spaces. This proves the claim.\\
	 
	Consider   the function $g_r(z)=\sum_{n\ge  0}\frac{ a_n}{n+1}(rz)^{n+1}$. In this case we have $g_r\in \dot H^{1/2}(S^1)$. We have
	$r \cdot f_r(z)=g'_r(z)$.
	 Since $g$ is harmonic in $\D$ we have 
	\begin{equation}\label{estL4}
	\|f_r\|_{L^2(\D)}\simeq\|\nabla g_r\|_{L^{2}{(\D)}}\lesssim \|g_r\|_{\dot H^{1/2}}\equiv \|f_r\|_{ H^{-1/2}}
	\end{equation}
	We let $r\to 1$ in \eqref{estL4} and get
	\begin{equation}\label{estL5}
	\|f\|_{L^2(\D)} \lesssim \| f \|_{H^{-1/2}(S^1)}. 
	\end{equation}
	Both inequalities  \eqref{embHA} and \eqref{embHA2} have been proved.


\begin{thebibliography}{9}
	
	\bibitem{bourgain1}
Bourgain, J. ; Brezis, H. \textit{ On the equation  $\div Y=f$ and application to control of phases. } J. Amer. Math. Soc. 16 (2003), no. 2, 393--426.
	
	\bibitem{bourgain2}
	Bourgain, J., Brezis, H.,
		\textit{New estimates for elliptic equations and Hodge type systems.} J. Eur. Math. Soc. (JEMS) 9 (2007), no. 2, 277--315. 
	
	\bibitem{Brez}
 Brezis, H. \textit{ Functional analysis, Sobolev spaces and partial differential equations.} Universitext. Springer, New York, 2011. xiv+599 pp.
	\bibitem{dalioriv}
	Da Lio, F., Rivi\`ere, T.,
	\textit{Critical Chirality in Elliptic Systems}, https://arxiv.org/abs/1907.10520.
	\bibitem{De} Delort, J.M. \textit{   Existence de nappes de tourbillon en dimension deux.} (French) [Existence of vortex sheets in dimension two] J. Amer. Math. Soc. 4 (1991), no. 3, 553-586.
\bibitem{DS} Duren, P; 
 Schuster, A. \textit{  Bergman Spaces.}  American Mathematical Soc., 2004.
	\bibitem{gilbert}
	Gilbert, J., Murray, M.,
	\textit{Clifford Algebras and Dirac Operators in Harmonic Analysis}. Cambridge Studies in Advanced Mathematics, 26. Cambridge University Press, Cambridge, 1991. viii+334.
	
	\bibitem{hamilton}Hamilton, M.,\textit{Mathematical Gauge Theory, 
with applications to the standard model of particle physics.} Universitext. Springer, Cham, 2017. xviii+657 
	\bibitem{mazya}  Maz'ya, V.  \textit{ Bourgain-Brezis type inequality with explicit constants. Interpolation theory and applications,} 247--252, Contemp. Math., 445, Amer. Math. Soc., Providence, RI, 2007.

	 	
	\bibitem{mironescu}
	Mironescu, P.,
	\textit{On some inequalities of Bourgain, Brezis, Maz'ya, and Shaposhnikova related to $ L^1$ vector fields. }
C. R. Math. Acad. Sci. Paris 348 (2010), no. 9--10.
\bibitem{Pel} Peloso, M. \textit{Classical spaces of holomorphic functions}  Appunti per il corso Argomenti
Avanzati di Analisi Complessa , Corso  di Laurea in Mathematica dell'Universit\`a di Milano
2014. Print.
\bibitem{Riesz} Riez, F. \textit{\"Uber die Randwerte einer analytischen Funktion.} Mathematische Zeitschrift, volume 18, (1923), 87--95.
 \bibitem{Riv3}  Rivi\`ere, T. \textit{ Sequences of smooth global isothermic immersions.} Comm. Partial Differential Equations 38 (2013), no. 2, 276-303.
\bibitem{ssv} Schikorra, A; Spector, D; Van Schaftingen, J. \textit{An L1-type estimate for Riesz potentials.} Rev. Mat. Iberoam. 33 (2017), no. 1, 291--303.
\bibitem{sw}   Stein E.M.;  Weiss G, \textit{On the theory of harmonic functions of several variables. I. The theory of $H^p$ spaces.} Acta Math. 103 (1960), 25--62.
	\bibitem{vs} Van Schaftingen, J. \textit{ Limiting Bourgain-Brezis estimates for systems of linear differential equations: theme and variations.} J. Fixed Point Theory Appl. 15 (2014), no. 2, 273--297.
 
 	\end{thebibliography}
\end{document}